%% file: main.tex
\documentclass[11pt]{article}
\usepackage{amsfonts,amsmath,amsthm,amssymb,mathdots}
\usepackage[all]{xy}
\usepackage[dvipsnames]{xcolor}
\usepackage[hidelinks,pdftex]{hyperref} 
\usepackage[bottom]{footmisc} 

\title{The Lichtenbaum--Quillen dimension \\ of complex varieties}
\author{Nicolas Addington and Elden Elmanto}
\date{}

\newcommand \A {\mathcal A}
\newcommand \C {\mathbb C}
\renewcommand \H {\mathcal H}
\renewcommand \P {\mathbb P}
\newcommand \Q {\mathbb Q}
\newcommand \Z {\mathbb Z}

\let \norwegianO \O
\DeclareRobustCommand \O {\ifmmode \mathcal O \else \norwegianO \fi}

\newcommand \dimQL {\dim_\mathrm{LQ}}
\newcommand \Dbcoh {D^b_\mathrm{coh}}

\DeclareMathOperator \Spec {Spec}
\DeclareMathOperator \Griff {Griff}
\DeclareMathOperator \Gr {Gr}
\DeclareMathOperator \Hom {Hom}
\DeclareMathOperator \Ext {Ext}
\DeclareMathOperator \Br {Br}
\DeclareMathOperator \Fil {Fil}
\DeclareMathOperator \gr {gr}
\DeclareMathOperator \RGamma {R\Gamma}
\DeclareMathOperator \Pic {Pic}
\DeclareMathOperator \NS {NS}
\DeclareMathOperator \rank {rank}
\DeclareMathOperator \Pf {Pf}
\DeclareMathOperator \Sym {Sym}
\DeclareMathOperator \cone {cone}
\DeclareMathOperator \Perf {Perf}

\newtheorem {thm} {Theorem}
\newtheorem {prop} {Proposition}
\numberwithin {prop} {section}
\newtheorem {cor} [prop] {Corollary}
\newtheorem {lem} [prop] {Lemma}
\theoremstyle{definition}
\newtheorem {defn} {Definition}

\newtheorem {examples} [prop] {Examples}

\begin{document}
\maketitle

\input intro
\input generalities
\input surfaces

\input 3folds
\input 4folds
\input 5folds

\bibliographystyle{plain}
\bibliography{ql}

\scriptsize
\noindent Nicolas Addington \\
Department of Mathematics \\
University of Oregon \\
Eugene, Oregon 97403-1222 \\
United States \\
adding@uoregon.edu \\

\noindent Elden Elmanto \\
Department of Mathematics \\
University of Toronto \\
40 St.\ George St. \\
Toronto, Ontario M5S 2E4 \\
Canada \\
elden.elmanto@utoronto.ca \\


\end{document}

%% file: intro.tex

\begin{abstract}
The Lichtenbaum--Quillen conjecture for smooth complex varieties states that algebraic and topological K-theory with finite coefficients become isomorphic in high degrees.  We define the ``Lichtenbaum--Quillen dimension'' of a variety in terms of the point where this happens, show that it is surprisingly computable, and analyze many examples.  It gives an obstruction to rationality, but one that turns out to be weaker than unramified cohomology and some related birational invariants defined by Colliot-Th\'el\`ene and Voisin using Bloch--Ogus theory.  Because it is compatible with semi-orthogonal decompositions, however, it allows us to prove some new cases of the integral Hodge conjecture using homological projective duality, and to compute the higher algebraic K-theory of the Kuznetsov components of the derived categories of some Fano varieties.
\end{abstract}

\section*{Introduction}

Let $X$ be a smooth complex projective variety, and consider the comparison map from algebraic to topological K-theory:
\[ \eta_n\colon K_n(X) \to KU^{-n}(X). \]
The image and kernel of $\eta_0$ are closely related to the image and kernel of the cycle class map $CH^*(X) \to H^{2*}(X,\Z)$, which are the subject of deep conjectures.  For $n \ge 1$, the image of $\eta_n$ is torsion \cite[\S6.3]{gillet}, so $\eta_n = 0$ in the many examples where there is no torsion in $H^*(X,\Z)$ and thus in $KU^*(X)$.  Moreover there is no Bott periodicity on $K_*(X)$.

With finite coefficients, the situation is completely different, as Thomason first showed in \cite{thomason_etale}.  To state his results, we introduce the Bott maps
\[ \tau_n\colon K_n(X,\Z/m) \to K_{n+2}(X,\Z/m), \]
roughly given by multiplication with a generator of $K_2(\Spec \C,\Z/m) = \Z/m$ but with subtleties that we review in \S\ref{bott_section}.  These are compatible with the comparison maps $\eta_n$ and the Bott isomorphisms on $KU^*$, given by multiplication with the corresponding generator $\beta \in KU^{-2}(\text{point}, \Z/m) = \Z/m$, as shown:
\[ \xymatrix{
K_n(X,\Z/m) \ar[d]_{\eta_n} \ar[r]^{\tau_n} & K_{n+2}(X,\Z/m) \ar[d]_{\eta_{n+2}} \\
KU^{-n}(X,\Z/m) \ar[r]^-{\cdot \beta}_-\cong & KU^{-n-2}(X,\Z/m).
} \]
Thomason proved that after taking direct limits over the Bott maps, the comparison maps become isomorphisms.  The Lichtenbaum--Quillen conjecture \cite[Thm.~4.1]{pw_varieties}, which is now a consequence of the Bloch--Kato conjecture or norm residue isomorphism theorem proved by Rost and Voevodsky \cite{hw}, gives more precise information: the comparison maps $\eta_n$ are isomorphisms for all $n \ge \dim(X) - 1$ and injective for $n = \dim(X) - 2$, or equivalently, the Bott maps $\tau_n$ are isomorphisms for all $n \ge \dim(X) - 1$ and injective for $n = \dim(X) - 2$.

But for some varieties, $\eta_n$ and $\tau_n$ become isomorphisms earlier than promised by the Lichtenbaum--Quillen conjecture.  For example, if $X \to Y$ is a $\P^r$ bundle (in the Zariski topology) then the K-theory of $X$ is just $r+1$ copies of the K-theory of $Y$, so the maps $\eta_n$ and $\tau_n$ are isomorphisms for $n \ge \dim(Y) - 1$ and injective for $n = \dim(Y) - 2$.  And we will see in \S\ref{semi-orth} that if $X$ is rational and $\dim(X) \ge 2$, then $\eta_n$ and $\tau_n$ are isomorphisms for $n \ge \dim(X) - 3$ and injective for $n = \dim(X) - 4$.

The goal of this paper is to explore this phenomenon and its interaction with birational geometry, algebraic cycles, and derived categories of coherent sheaves.  We begin by introducing the following definition.
\begin{defn} \label{QL_dim}
The \emph{Lichtenbaum--Quillen dimension} of a smooth complex projective variety $X$, denoted $\dimQL(X)$, is the smallest integer $d$ such that for all $m$, the comparison maps
\[ \eta_n\colon K_n(X,\Z/m) \to KU^{-n}(X,\Z/m) \]
are isomorphisms for $n \ge d-1$ and injective for $n = d-2$, or equivalently, the Bott maps
\[ \tau_n\colon K_n(X,\Z/m) \to K_{n+2}(X,\Z/m) \]
are isomorphisms for $n \ge d-1$ and injective for $n = d-2$.
\end{defn}
Thus the statements above amount to saying that $\dimQL(X) \le \dim(X)$ for all $X$, but if $X \to Y$ is a $\P^r$ bundle then $\dimQL(X) = \dimQL(Y)$, and if $X$ is rational and $\dim(X) \ge 2$ then $\dimQL(X) \le \dim(X) - 2$.

It is easy to see that when $X$ is a curve, we have $\dimQL(X) = 0$ if and only if the genus is 0: we always have $K_{-1}(X,\Z/m) = 0$, so the comparison map $\eta_{-1}$ cannot be an isomorphism unless $KU^1(X,\Z/m) = H^1(X,\Z/m)$ vanishes.  When $\dim(X) \ge 2$, it looks at first glance as if $\dimQL(X)$ might be totally uncomputable, but we introduce a spectral sequence, related to the motivic spectral sequence, that computes the cone of the Bott maps $\tau_n$.\footnote{It is interesting to note that the cofiber of $\tau$ has another life in motivic homotopy theory \cite{gwx}, where it has been used to compute homotopy groups of spheres \cite{iwx}.}  \begin{thm} \label{spec_seq}
There is a spectral sequence
\[ E_2^{p,q} = H^p(\H^{-q}_X(\Z/m)) \Longrightarrow (K/\tau)_{-p-q}(X,\Z/m), \]
where $H^*(\H^*_X(\Z/m))$ are the Bloch--Ogus cohomology groups discussed in \S\ref{Hi_Hj_sec} and \S\ref{spec_seq_surfaces}, and the groups $(K/\tau)_*(X,\Z/m)$ are defined carefully in \S\ref{bott_section} and fit into the long exact sequence
\begin{multline*}
\dotsb \to K_{n-2}(X,\Z/m) \xrightarrow{\tau_{n-2}} K_n(X,\Z/m) \to (K/\tau)_n(X,\Z/m) \\
\to K_{n-3}(X,\Z/m) \xrightarrow{\tau_{n-3}} K_{n-1}(X,\Z/m) \to \dotsb.
\end{multline*}
\end{thm}
\noindent Our $E_2$ page has the same terms as the $E_2$ page of the Bloch--Ogus--Leray spectral sequence, which abuts to the coniveau filtration on $H^*_\text{sing}(X,\Z/m)$, but with different indexing or (equivalently) with the differentials going a different way.  This is reminiscent of other spectral sequences that run one way to give a geometric invariant and another way to give a derived invariant, like the Hodge--de\,Rham spectral sequence and the Hochschild--Kostant--Rosenberg spectral sequence, which compute de Rham cohomology and Hochschild homology from $H^*(\Omega^*_X)$, or in positive characteristic, the slope spectral sequence and the descent spectral sequence, which compute crystalline cohomology and Hesselholt's TR invariants from Hodge--Witt cohomology \cite{antieau-bragg}.

\subsection*{Geometric conditions and rationality}

Because of this close connection with Bloch--Ogus groups, our computations of Lichtenbaum--Quillen dimension end up involving unramified cohomology and some related birational invariants studied by Colliot-Th\'el\`ene and Voisin in \cite[Props.~3.3 and 3.4]{ctv}.  We obtain a number of concrete statements about Lichtenbaum--Quillen dimension for surfaces, rationally connected 3-folds and 4-folds, and cubic 5-folds.

\begin{thm} \label{surface_thm}
Let $X$ be a smooth complex projective surface.
\begin{enumerate}
\item $\dimQL(X) \le 1$ if and only if the geometric genus $p_g(X) = h^{0,2}(X)$ vanishes and $H^3(X,\Z)$ is torsion-free.
\item $\dimQL(X) = 0$ if and only if in addition, the irregularity $q(X) = h^{0,1}(X)$ vanishes.
\end{enumerate}
\end{thm}
\noindent Thus while we said above that a rational surface has $\dimQL(X) = 0$, we see here that the converse is not true: Barlow surfaces \cite{barlow} and Dolgachev surfaces \cite{dolgachev} have $p_g = q = 0$ and are simply connected, but are not rational.\bigskip

In general we will see that if $X$ is a $d$-dimensional variety with $h^{0,d}(X) \ne 0$, then $\dimQL(X) = d$.  But starting in dimension 3 we will restrict our attention to rationally connected varieties.  Recall that a smooth complex projective variety is \emph{rationally connected} if any two points can be joined by a chain of rational curves \cite[Cor.~4.28]{debarre_book}.  This includes all unirational varieties, as well as all smooth Fano varieties \cite{kmm}, and thus all $d$-dimensional smooth hypersurfaces of degree $\le d$.  It implies that $CH_0(X) = \Z$, but is strictly stronger. 
It implies that $h^{0,i}(X) = 0$ for all $i > 0$ \cite[Cor.~4.18(a)]{debarre_book}.

\begin{thm} \label{3fold_thm}
Let $X$ be a smooth complex projective 3-fold that is rationally connected.
\begin{enumerate}
\item $\dimQL(X) \le 2$
\item $\dimQL(X) \le 1$ if and only if $H^3(X,\Z)$ is torsion-free.
\item $\dimQL(X) = 0$ if and only if $H^3(X,\Z) = 0$.
\end{enumerate}
\end{thm}
\noindent Again we see that the condition $\dimQL(X) \le 1$ does not give a very strong obstruction to rationality: all smooth Fano 3-folds have torsion-free cohomology, even though many of them are irrational, including cubic and quartic 3-folds.  The simplest 3-fold for which Lichtenbaum--Quillen dimension \emph{does} obstruct rationality is the Artin--Mumford 3-fold \cite{am}, although this does not amount to a new proof.
\pagebreak 

\begin{thm} \label{4fold_thm}
Let $X$ be a smooth complex projective 4-fold that is rationally connected.
\begin{enumerate}
\item $\dimQL(X) \le 3$ if and only if algebraic and homological equivalence coincide on $CH_1(X)$, and in the coniveau filtration\footnote{Recall that $N^c H^l(X,A)$ is defined as the set of classes $\alpha \in H^l(X,A)$ for which there is a Zariski closed set $Z \subset X$ of codimension $\ge c$ with $\alpha|_{X\setminus Z} = 0$.
If $c > l/2$ then $N^c = 0$,
and if $c \le l - \dim(X)$ then $N^c$ is all of $H^l(X,A)$.
If $l = 2c$ then $N^c H^{2c}(X,\Z)$ is the set of algebraic classes.}
\[ N^2 H^5(X,\Z) \subset N^1 H^5(X,\Z) = H^5(X,\Z), \]
the inclusion is an equality.

\item $\dimQL(X) \le 2$ if and only if in addition, the integral Hodge conjecture holds for $H^4(X,\Z)$ and $H^6(X,\Z)$.

\item If in addition $h^{1,3}(X) = 0$, and $H^5(X,\Z)$ and $H^3(X,\Z)$ are torsion-free, then $\dimQL(X) \le 1$.

\item[(c$'$)] If $\dimQL(X) \le 1$ then $h^{1,3}(X) = 0$, and $H^5(X,\Z)$ is torsion-free.

\end{enumerate}
\end{thm}

\noindent We do not know whether $\dimQL(X) \le 1$ implies that $H^3(X,\Z)$ is torsion-free.  A characterization of rationally connected 4-folds with $\dimQL(X) = 0$ would become too involved.

The invariants appearing in Theorem \ref{4fold_thm}(a) and (b) are more or less familiar birational invariants.  The integral Hodge conjecture on $H^4(X,\Z)$ and $H^{2d-2}(X,\Z)$ has been studied in the context of birational geometry by Voisin and her coauthors in \cite[Lem.~1]{voisin-soule}, \cite[Lem.~15]{voisin_Z_hodge}, and especially \cite{ctv}; note that for any Fano 4-fold, $H^6(X,\Z)$ is generated by algebraic classes by \cite[Thm.~1.7(ii)]{voisin-horing}.  The Griffiths group $\Griff_1(X) := CH_1(X)_\text{hom}/CH_1(X)_\text{alg}$ is discussed in \cite[\S2.1]{voisin_lueroth} and \cite[\S1.3.2]{voisin_birational}, and the quotient $H^5/N^2 H^5$ figures prominently in \cite{voisin_deg_4}.

We will apply Theorem \ref{4fold_thm} to the two families of Fano 4-folds shown by Hassett, Pirutka, and Tschinkel to contain both rational and irrational members, and will see that all members of both families have $\dimQL(X) = 2$.  We will see that all cubic 4-folds and Gushel--Mukai 4-folds have $\dimQL(X) = 2$, so Lichtenbaum--Quillen dimension gives no obstruction to their rationality.  We will see that quartic 4-folds that are very general in the Noether--Lefschetz sense have $\dimQL(X) = 2$; Totaro \cite{totaro_irr} proved that a very general quartic 4-fold is irrational, but we do not know how the latter ``very general'' interacts with the former.

On the other hand, we will see that the unirational 4-fold constructed by Schreieder in \cite[Cor.~1.6]{schreieder_small_slopes} has $\dimQL(X) = 3$. \bigskip

In dimension 5, we will analyze one example, relying on calculations of Fu and Tian \cite{fu_tian}:
\begin{thm} \label{cubic_5fold}
If $X$ is a smooth complex cubic 5-fold, then $\dimQL(X) = 1$.
\end{thm}

\subsection*{Lichtenbaum--Quillen dimension as a derived invariant}

The algebraic and topological K-theory of $X$ depend only on the derived category of coherent sheaves $\Dbcoh(X) = \Perf(X)$, as we review in \S\ref{semi-orth}, so the same is true of $\dimQL(X)$.  In fact one can extend Definition \ref{QL_dim} to any $\C$-linear dg-category or stable $\infty$-category, using Blanc's topological K-theory \cite{blanc}.  When the derived category of a variety admits a semi-orthogonal decomposition
\[ \Dbcoh(X) = \langle \A_1, \dotsc, \A_k \rangle, \]
we will see that
\[ \dimQL(X) = \max(\dimQL(\A_1),\, \dotsc,\, \dimQL(\A_k)). \]
We extract two kinds of results from this observation. \bigskip

First, we use homological projective duality to produce some interesting Fano 4-folds with $\dimQL(X) = 2$, which therefore satisfy all the conditions of Theorem \ref{4fold_thm}(a) and (b):

\begin{thm} \label{hpd_thm}
The following Fano 4-folds $X$ satisfy the integral Hodge conjecture on $H^4(X,\Z)$, and $\Griff_1(X) = 0$.
\begin{enumerate}
\item Intersections of $\Gr(2,7)$ with 6 general hyperplanes in the Pl\"ucker embedding.
\item ``Pfaffian'' 4-folds obtained as linear sections of the space of $7 \times 7$ skew-symmetric matrices of rank 4.
\item The linear sections of the double quintic symmetroid studied by Ottem and Rennemo in \cite{ottem-rennemo}.
\end{enumerate}
\end{thm}
\noindent The third example has $H^3(X,\Z) = \Z/2$ and thus is not even stably rational.  Rationality of the first two examples does not seem to have been studied.\footnote{Shortly before this paper was accepted, we learned that the first example was in fact known to be rational: according to Russo and Staglian\`o \cite[Table 1, line ix]{russo-stagliano}, it is birational to a Pfaffian cubic fourfold.  We thank N.A.'s student Zhe Xu for bringing this reference to our attention.  We still do not know about rationality of the second example.}

In each case, $\Dbcoh(X)$ admits a semi-orthogonal decomposition consisting of an exceptional collection and the derived category of a surface $S$ with $h^{0,2} \ne 0$.  In the first example, where $H^*(X,\Z)$ is torsion-free, it is easier to deduce the integral Hodge conjecture using Perry's \cite[Prop.~5.16(2)]{perry}, but the result on $\Griff_1(X)$ is new.  In the second example, we do not know whether there is torsion in $H^*(X,\Z)$.  In all three examples, the fully faithful functors $\Dbcoh(S) \to \Dbcoh(X)$ are induced by ideal sheaves of explicit subschemes of $S \times X$, and it should be possible to prove Theorem \ref{hpd_thm} more directly by analyzing the geometry of those correspondences; but it is also gratifying to get integral information about $H^4$ and $CH_1$ from derived categories, not just rational information. \bigskip

Second, we compute the higher K-theory of some Kuznetsov components of derived categories.  For many Fano varieties $X$, Kuznetsov has identified a semi-orthogonal decomposition of $\Dbcoh(X)$ consisting of an exceptional collection and an interesting piece $\A_X$, which often behaves like the derived category of a lower-dimensional variety.  Most famously, if $X$ is a cubic 4-fold \cite{kuz_cubics} or Gushel--Mukai 4-fold \cite{kuz_perry}, then the Kuznetsov component $\A_X$ behaves like the derived category of a K3 surface in terms of Serre duality, Hochschild homology and cohomology, and the Euler pairing on $KU^0$.  We find that its algebraic K-theory behaves like that of a K3 surface as well:
\begin{thm} \label{4fold_KAX}
Let $X$ be a smooth complex cubic 4-fold or Gushel--Mukai 4-fold, and let $\mathcal A_X \subset \Dbcoh(X)$ be Kuznetsov's K3 category.
\begin{enumerate}
\item $\dimQL(\A_X) = 2$, meaning that the Bott maps
\[ \tau_n\colon K_n(\A_X, \Z/m) \to K_{n+2}(\A_X, \Z/m) \]
are isomorphisms for $n \ge 1$ and injective for $n=0$.
\item With integral coefficients, we have
\[ K_n(\A_X) = \begin{cases}
\Z^\rho \oplus V_0 & \text{if $n=0$,} \\
(\Q/\Z)^{24} \oplus V_n & \text{if $n$ is odd and $n \ge 1$, and} \\
V_n & \text{if $n$ is even and $n \ge 2$,}
\end{cases} \]
where $\rho$ is the rank of $K_\mathrm{num}(\A_X)$ and the $V_n$s are (possibly infinite-dimensional) $\Q$-vector spaces.  In fact $V_0 = CH_1(X)_\mathrm{hom}$.
\end{enumerate}
\end{thm}
\noindent The algebraic K-theory of a K3 surface $S$ admits a similar description, with $V_0 = CH_0(S)_\text{hom}$, as we review in \S\ref{cubic_gm_sec}.  We wonder whether a similar conclusion holds for the K3 categories obtained from Debarre--Voisin Fano 20-folds -- that is, hyperplane sections of $\Gr(3,10)$ -- but these seem far out of reach.  The fact that Gushel--Mukai 4-folds have no torsion in $CH_1$ is new and of independent interest; the same proof applies to quartic 4-folds and several other Fano 4-folds.

We will also describe the K-theory of Kuznetsov components of derived categories of Fano 3-folds in \S\ref{fano_3folds}, and of cubic 5-folds in \S\ref{cubic_5fold_sec}.

\subsection*{Acknowledgements}

We thank Dan Dugger and Ben Antieau for extensive discussions.  Claire Voisin was generous with her advice, as were Jean-Louis Colliot-Th\'el\`ene, Lie Fu, James Hotchkiss, John Christian Ottem, J\o rgen Rennemo, Stefan Schreieder, and Zhiyu Tian.  N.A.\ thanks the Max Planck Institute for Mathematics in Bonn for their hospitality, and was supported by NSF grant no.\ DMS-1902213.  E.E. was supported by the NSERC grant RGPIN-2025-07114, ``Motivic cohomology: theory and applications.''  We thank the referees for their careful reading and many helpful suggestions.

\subsection*{Conventions}

For spectra, we use the notation and syntax of higher algebra \cite{lurie-ha}; we write $\otimes$ and $\oplus$ and $[1]$ where some authors would write $\wedge$ and $\vee$ and $\Sigma$.  Algebraic K-theory for us is \emph{connective} algebraic K-theory, that is, the universal additive invariant in the sense of \cite{bgt}.

%% file: generalities.tex

\section{General Results}

\subsection{The Bott map and its cofiber} \label{bott_section}

Let $X$ be a variety over the complex numbers.  In this section we explain the Bott maps 
\begin{equation} \label{tau_n}
\tau_n\colon K_n(X, \Z/m) \to K_{n+2}(X, \Z/m),
\end{equation}
their relation to the comparison maps 
\begin{equation} \label{eta_n}
\eta_n\colon K_n(X, \Z/m) \to KU^{-n}(X, \Z/m),
\end{equation}
and the equivalence between formulating Definition \ref{QL_dim} with \eqref{tau_n} or \eqref{eta_n}.  Then we upgrade the Bott maps \eqref{tau_n} to a map of spectra, take its cofiber, and give another equivalent formulation of Definition \ref{QL_dim}. \bigskip

Topological K-theory satisfies Bott periodicity: if $\beta$ is a generator of $KU^{-2}(\text{point}) \cong \Z$, then multiplication by $\beta$ gives isomorphisms
\[ \xymatrix{
KU^{-n}(X) \ar[r]^-{\cdot\beta}_-\cong & KU^{-n-2}(X)
} \]
and
\[ \xymatrix{
KU^{-n}(X,\Z/m) \ar[r]^-{\cdot\beta}_-\cong & KU^{-n-2}(X,\Z/m)
} \]
for all $n$ and $m$.

In algebraic K-theory with $\Z$ coefficients, there is no Bott periodicity, but with $\Z/m$ coefficients we can proceed as follows.  Suslin's rigidity theorem \cite[Cor.~4.7]{suslin_local_fields} implies that $\eta_n\colon K_n(\Spec \C,\Z/m) \to KU^{-n}(\text{point},\Z/m)$ is an isomorphism for all $n \ge 0$.  We would like to take the same generator $\beta \in KU^{-2}(\text{point})$, map it to a generator of $KU^{-2}(\text{point}, \Z/m) \cong \Z/m$, lift it to $K_2(\Spec \C,\Z/m)$, and let it act on $K_*(X,\Z/m)$ to give the maps \eqref{tau_n} above.

But there are well-known difficulties when $m \equiv 2 \pmod 4$, which we cannot ignore because $\Z/2$ coefficients are the most interesting ones for our applications.  In that case the Moore spectrum $\mathbb S/m$ does not admit a unital multiplication $\mathbb S/m \otimes \mathbb S/m \to \mathbb S/m$, so $K(\Spec \C,\Z/m)$ and $K(X,\Z/m)$, which are obtained from $K(\Spec \C)$ and $K(X)$ by tensoring with $\mathbb{S}/m$, may also fail to be ring spectra.\footnote{Pedrini and Weibel \cite{pw_varieties} cite Araki and Toda \cite[Thm.~10.7]{at2} to say that the homotopy groups $K_*(\Spec \C,\Z/m)$ and $K_*(X,\Z/m)$ are nonetheless graded rings because $\sqrt{-1} \in \C$, but soon we will really need to work with spectra.}  But the action $\mathbb S/2m \otimes \mathbb S/m \to \mathbb S/m$ is better behaved, so we define  \eqref{tau_n} by letting the image of $\beta$ in $K_2(\Spec \C,\Z/2m)$ act on $K_*(X,\Z/m)$.  There are further difficulties with associativity of this action when $m$ is divisible by 4 but not 16, or by 3 but not 9, but these do not affect us.  For further discussion of these issues, see Thomason's \cite[App.~A]{thomason_etale}, or for a more recent perspective see \cite[\S3 and App.~A]{elden_bott} or \cite{burklund}.

Thus we have defined the Bott maps \eqref{tau_n} for all $m$, and we have a commutative square \[ \xymatrix{
K_n(X,\Z/m) \ar[r]^-{\tau_n} \ar[d]_{\eta_n} & K_{n+2}(X,\Z/m) \ar[d]_{\eta_{n+2}} \\
KU^{-n}(X,\Z/m) \ar[r]^-{\cdot \beta}_-\cong & KU^{-n-2}(X,\Z/m).
} \]
Knowing that $\eta_n$ is an isomorphism for $n \gg 0$, we see for any given $N$ that $\eta_n$ is an isomorphism (or injective) for all $n \ge N$ if and only if $\tau_n$ is an isomorphism (or injective) for all $n \ge N$: thus the two formulations of Definition \ref{QL_dim} are equivalent, as promised. \bigskip

Next we upgrade the Bott maps \eqref{tau_n} to a map of spectra
\begin{equation} \label{bott_on_spectra}
\tau\colon  K(X,\Z/m)[2] \to K(X,\Z/m)
\end{equation}
When $m \equiv 2 \pmod 4$, we start with our generator of $K_2(\Spec \C,\Z/2m)$ and take the corresponding map
\[ S^2 \to K(\Spec \C,\Z/2m). \]
Then we tensor with $K(X,\Z/m)$ to get
\[ S^2 \otimes K(X,\Z/m) \to K(\Spec \C,\Z/2m) \otimes K(X,\Z/m), \]
and post-compose with the map
\[ K(\Spec \C,\Z/2m) \otimes K(X,\Z/m) \to K(X,\Z/m) \]
obtained by tensoring the multiplication map
\[ K(\Spec \C) \otimes K(X) \to K(X) \]
with the action
\[ \mathbb{S}/2m \otimes \mathbb{S}/m \to \mathbb{S}/m. \]
The case $m \not\equiv 2 \pmod 4$ is more straightforward: we do the same with $m$ in place of $2m$.  

Finally we define $K/\tau(X,\Z/m)$ as the cofiber of the spectrum-level Bott map \eqref{bott_on_spectra}.  Its homotopy groups $(K/\tau)_*(X,\Z/m)$ fit into a long exact sequence
\begin{multline*}
\dotsb \to K_{n-2}(X,\Z/m) \xrightarrow{\tau_{n-2}} K_n(X,\Z/m) \to (K/\tau)_n(X,\Z/m) \\
\to K_{n-3}(X,\Z/m) \xrightarrow{\tau_{n-3}} K_{n-1}(X,\Z/m) \to \dotsb,
\end{multline*}
giving another equivalent formulation of Definition \ref{QL_dim}:
\begin{defn} \label{QL_dim'}
The \emph{Lichtenbaum--Quillen dimension} of a smooth complex projective variety $X$ is the smallest integer $d$ such that for all $m$, the cofiber of the spectrum-level Bott map satisfies
\[ (K/\tau)_n(X,\Z/m) = 0 \]
for $n > d$.
\end{defn}

\subsection{Derived categories and semi-orthogonal decompositions} \label{semi-orth}

If $X$ and $Y$ are smooth complex projective varieties then every equivalence $\Dbcoh(X) \to \Dbcoh(Y)$ is induced by a Fourier--Mukai kernel, which induces weak equivalences $K(X) \to K(Y)$ and $KU(X) \to KU(Y)$, and similarly with finite coefficients, compatible with all comparison and Bott maps.  Thus $\dimQL(X)$ depends only on $\Dbcoh(X)$.

In fact, if $\mathcal C$ is any dg-category or stable $\infty$-category over $\mathbb C$, we can define the spectrum $K(\mathcal C)$ using Waldhausen's S-construction as in \cite{bgt}, and $KU(\mathcal C)$ following Blanc \cite{blanc}, and with finite coefficients we can define Bott maps just as in the previous section, so we can extend Definition \ref{QL_dim}:
\begin{defn}
The \emph{Lichtenbaum--Quillen dimension} of a small idempotent-complete $\C$-linear pretriangulated dg-category $\mathcal C$, or equivalently a small idempotent-complete $\C$-linear stable $\infty$-category $\mathcal C$, is the smallest integer $d$ such that any of the following equivalent conditions holds for all $m$:
\begin{enumerate}
\item The comparison maps $\eta_n\colon K_n(\mathcal C,\Z/m) \to KU^{-n}(\mathcal C,\Z/m)$ are isomorphisms for $n \ge d-1$ and injective for $n = d-2$.
\item The Bott maps $\tau_n\colon K_n(\mathcal C,\Z/m) \to K_{n+2}(\mathcal C,\Z/m)$ are isomorphisms for $n \ge d-1$ and injective for $n = d-2$.
\item The cofiber of the spectrum-level Bott map satisfies $(K/\tau)_n(\mathcal C, \Z/m) = 0$ for $n > d$.
\end{enumerate}
If no such $d$ exists, we define $\dimQL(\mathcal C) = \infty$.
\end{defn}
The equivalence between (a) and (b) depends on a result of Antieau and Heller \cite[Thm.~3.3]{antieau-heller}: they proved Blanc's conjecture \cite[Conj.~4.27]{blanc} that the natural map from (connective) algebraic K-theory with finite coefficients to semitopological K-theory with finite coefficients is an equivalence, and this implies that $KU(\mathcal C, \Z/m)$ is the direct limit of the Bott maps on $K(\mathcal C, \Z/m)$.  When $\mathcal C = \Dbcoh(X) = \Perf(X)$ we recover our earlier definitions.

In general we do not know whether $\dimQL(\mathcal C)$ is finite, even when $\mathcal C$ is smooth and proper; cf.\ \cite[Question~4.2]{antieau-heller}.  But in practice we are only interested in admissible subcategories $\mathcal C \subset \Dbcoh(X)$, in which case we will see that $\dimQL(\mathcal C) \le \dimQL(X)$.  Recall that a full triangulated subcategory of $\Dbcoh(X)$ is called \emph{admissible} if the inclusion functor has left and right adjoints; it inherits its dg or stable $\infty$ enhancement from $\Dbcoh(X)$.

A \emph{semi-orthogonal decomposition}
\[ \Dbcoh(X) = \langle \A_1, \dotsc, \A_k \rangle \]
is a sequence of admissible subcategories, with no Homs or Exts from right to left, that generate $\Dbcoh(X)$ as a triangulated category.
\begin{examples} \label{sod_examples} \ 
\begin{enumerate}
\item
Beilinson's exceptional collection
\[ \Dbcoh(\P^r) = \langle \O_{\P^r}, \O_{\P^r}(1), \dotsc, \O_{\P^r}(r) \rangle. \]
Here each $\O_{\P^r}(i)$ is a shorthand for the subcategory that it generates, which is an admissible subcategory equivalent to $\Dbcoh$(point) because $\O_{\P^r}(i)$ is an \emph{exceptional object}, with $\Hom(\O_{\P^r}(i), \O_{\P^r}(i)) = \C$ and $\Ext^j(\O_{\P^r}(i), \O_{\P^r}(i)) = 0$ for $j \ne 0$.
\item More generally, if $P \to X$ is a $\P^r$ bundle in the Zariski topology, Orlov has given a semi-orthogonal decomposition
\[ \Dbcoh(P) = \big\langle \overbrace{\Dbcoh(X),\, \dotsc,\, \Dbcoh(X)}^\text{$r+1$ times} \big\rangle. \]
The notation means that we have $r+1$ fully faithful functors \linebreak 
$\Dbcoh(X) \to \Dbcoh(P)$, and their images give a semi-orthogonal decomposition.  A semi-orthogonal decomposition induces a direct sum decomposition on $K(X)$, so this refines Quillen's \cite[\S8 Thm.~2.1]{quillen}.
\item If $P \to X$ is only a $\P^r$ bundle in the \'etale or analytic topology, Bernardara \cite{bernardara} has given a semi-orthogonal decomposition
\[ \Dbcoh(P) = \langle \Dbcoh(X), \Dbcoh(X,\alpha), \Dbcoh(X,2\alpha), \dotsc, \Dbcoh(X,r\alpha) \rangle, \]
which refines Quillen's \cite[\S8 Thm.~4.1]{quillen}; here $\alpha \in \Br(X)$ is the Brauer class obtained as the obstruction to writing $P$ as the projectivization of a vector bundle, and $\Dbcoh(X,\alpha)$ is the derived category of $\alpha$-twisted sheaves.
\item If $\tilde X$ is the blow-up of $X$ along a smooth center $Z$, Orlov has given a semi-orthogonal decomposition
\[ \Dbcoh(\tilde X) = \langle \Dbcoh(X), \overbrace{\Dbcoh(Z), \dotsc, \Dbcoh(Z)}^\text{codim\,$Z - 1$ times} \rangle, \]
which refines Thomason's \cite[Thm.~2.1]{thomason_blowup}.
\item If $X$ is a smooth complete intersection of two quadrics in $\P^{2g+1}$, Bondal and Orlov \cite[Thm.~2.9]{bo} have given a semi-orthogonal decomposition
\[ \Dbcoh(X) = \langle \Dbcoh(C), \O_X, \O_X(1), \dotsc, \O_X(2g-1) \rangle, \]
which refines Reid's \cite[Thm.~4.14]{reid}; here $C$ is a hyperelliptic curve of genus $g$ whose construction goes back to Weil.  We will encounter many semi-orthogonal decompositions of this kind, consisting of one ``interesting piece'' and an exceptional collection that does not vary when $X$ deforms.
\end{enumerate}
\end{examples}

A semi-orthogonal decomposition on $\Dbcoh(X)$ induces a canonical direct sum decomposition on $K(X)$ and $KU(X)$: one approach is to take the projection kernels of \cite[Thm.~7.1]{kuz_basechange} and let them induce idempotents on $K(X)$ and $KU(X)$.  In fact one could take this as the definition of $KU$ of an admissible subcategory rather than invoking Blanc's machinery.  The idempotents and decompositions are compatible with all comparison and Bott maps, so we also get a direct sum decomposition of $K/\tau(X)$, and we find:
\begin{prop} \label{sod_vs_QL}
If $\Dbcoh(X) = \langle \A_1, \dotsc, \A_k \rangle$, then 
\begin{center} \phantom{\qed} \hfill
$\dimQL(X) = \max(\dimQL(\A_1),\, \dotsc,\, \dimQL(\A_k))$.
\hfill \qed \end{center}
\end{prop}
\begin{cor} \label{Pn_is_QL0}
$\dimQL(\P^r) = 0$.  More generally, if $P \to X$ is a $\P^r$ bundle in the Zariski topology, then $\dimQL(P) = \dimQL(X)$.  If it is only a $\P^r$ bundle in the \'etale topology, then $\dimQL(P) \ge \dimQL(X)$.
\end{cor}
\begin{proof}
This follows from Examples \ref{sod_examples}(a)--(c) and Proposition \ref{sod_vs_QL}.
\end{proof}
\noindent It is tempting to conjecture that equality holds for \'etale $\P^r$ bundles as well; we will return to this question in \S\ref{conic_section}.

\begin{cor} \label{K/tau_birational} \ 
\begin{enumerate}
\item For a smooth complex projective variety $X$ of dimension $d$, the groups $(K/\tau)_d(X,\Z/m)$ and $(K/\tau)_{d-1}(X,\Z/m)$ are birational invariants.
\item If $Y$ is birational to $X$ and either $\dimQL(X) = d$, or $\dimQL(X) = d - 1$, or $\dimQL(X) \le d - 2$, then the same is true of $Y$.
\item In particular, if $X$ is rational then $\dimQL(X) \le \max(d - 2,\,0)$.
\end{enumerate}
\end{cor}
\begin{proof}
(a) By weak factorization \cite{akmw}, every birational map between smooth complex projective varieties factors into a sequence of blow-ups and blow-downs along smooth centers, so it is enough to show that if $\tilde X$ is the blow-up of $X$ along a smooth center $Z$ of codimension $\ge 2$, then
\[ (K/\tau)_n(\tilde X,\Z/m) = (K/\tau)_n(X, \Z/m) \]
for $n=d$ and $d-1$.  This follows from the semi-orthogonal decomposition in Example \ref{sod_examples}(d) and the fact that $\dimQL(Z) \le \dim(Z) \le d-2$.

(b) This is immediate from (a) and Definition \ref{QL_dim'}.

(c) We have $\dimQL(\P^d) = 0$ by Corollary \ref{Pn_is_QL0}.
\end{proof}
\noindent But we remark that $\dimQL(X)$ is not itself a birational invariant: if we blow up $\P^3$ along a curve of genus $\ge 1$, then $\dimQL$ jumps from 0 to 1.

\begin{cor}
If $X$ is a smooth complete intersection of two quadrics in $\P^{2g+1}$, then $\dimQL(X) = 1$.
\end{cor}
\begin{proof}
This follows from Example \ref{sod_examples}(e) and Proposition \ref{sod_vs_QL}.
\end{proof}
This is as far as we can get with semi-orthogonal decompositions alone; to proceed further, we  construct a motivic spectral sequence for $K/\tau$, which allows us to make more subtle connections between Lichtenbaum--Quillen dimension and geometry.

\subsection{A motivic spectral sequence for \texorpdfstring{$K/\tau$}{K/tau}} \label{Hi_Hj_sec}

In this section we introduce a spectral sequence to compute $(K/\tau)_*(X,\Z/m)$.  To describe the $E_2$ page, let
\[ \pi\colon X_\text{an} \to X_\text{Zar} \]
be the identity map from $X$ in the analytic topology to the Zariski topology, and for any Abelian group $A$, let
\[ \H^i_X(A) := R^i \pi_*\,\underline A, \]
where $\underline A$ is the constant sheaf on $X$ in the analytic topology; in other words, $\H^i_X(A)$ is the sheaf associated to the presheaf which maps a Zariski open set $U \subset X$ to $H^i_\text{sing}(U,A)$.
{ \renewcommand{\thethm}{\ref{spec_seq}}
\addtocounter{thm}{-1}
\begin{thm}
There is a spectral sequence
\[ E_2^{p,q} = H^p(\H^{-q}_X(\Z/m)) \Longrightarrow (K/\tau)_{-p-q}(X,\Z/m), \]
functorial in $X$.
\end{thm}
}
\noindent We postpone the proof to the next section, but first discuss some consequences.

Bloch and Ogus \cite[eq.~(0.3)]{bloch-ogus} proved that $H^i(\H^j_X(A))$ vanishes for $i > j$ and any $A$, so we can write out the $E_2$ page of the spectral sequence from Theorem \ref{spec_seq} as follows, suppressing the $\Z/m$ coefficients:\footnote{It would be interesting to understand the differentials, which should be some incarnation of Steenrod squares or powers.  We suspect that many of them vanish.}
\begin{equation} \label{our_E_2}
\begin{split}
\xymatrix@=1ex{
\ar@{-}[]+<1.6em,-.8em>;[rrrrrr]+<0em,-.8em> 
\ar@{-}[]+<1.6em,-.8em>;[dddddd]+<1.6em,0em> 
& p=0 & 1 & 2 & 3 & 4 & \\
q=0 & H^0(\H^0) \\
-1 & H^0(\H^1) \ar[rrd] & H^1(\H^1) \\
-2 & H^0(\H^2) \ar[rrd] & H^1(\H^2) \ar[rrd] & H^2(\H^2) \\
-3 & H^0(\H^3) \ar[rrd] & H^1(\H^3) \ar[rrd] & H^2(\H^3) \ar[rrd] & H^3(\H^3) \\
-4 & H^0(\H^4) \ar[rrd] & H^1(\H^4) \ar[rrd] & H^2(\H^4) \ar[rrd] & H^3(\H^4) \ar[rrd] & H^4(\H^4) \\
& \vdots & \vdots & \vdots & \vdots & \vdots & \ddots
}
\end{split}
\end{equation}
If $\dim(X) = d$ then we also have $\H^j_X(A) = 0$ for $j > d$: the sheafification can be performed using \emph{affine} open sets, which are homotopy equivalent to CW complexes of (real) dimension at most $d$; this is even true for singular varieties \cite[Thm.~2.12]{karchyauskas} \cite{hamm}.
From this we see again that
\[ (K/\tau)_n(X,\Z/m) = 0 \qquad \text{for }n > d, \]
which is exactly the statement of the Lichtenbaum--Quillen conjecture for $X$; but we note that the proof of Theorem~\ref{spec_seq} will rely on Rost and Voevodsky's resolution of the Bloch--Kato conjecture.

We also get
\begin{equation} \label{K/tau_d}
(K/\tau)_d(X, \Z/m) = H^0(\H^d_X(\Z/m))
\end{equation}
and an exact sequence of low-degree terms
\begin{multline} \label{K/tau_d-1}
0 \to H^1(\H^d_X(\Z/m)) \to (K/\tau)_{d-1}(X, \Z/m) \\
\to H^0(\H^{d-1}_X(\Z/m)) \to H^2(\H^d_X(\Z/m))
\end{multline}
which we will use often. \medskip

The groups $H^0(\H^j_X(A))$ are sometimes called unramified cohomology and denoted $H^j_\text{nr}(X,A)$; they are well-known birational invariants, and if $X = \P^d$ then they vanish for $j > 0$.  Colliot-Th\'el\`ene and Voisin have shown that $H^i(\H^d_X(A))$ are birational invariants as well \cite[Prop.~3.4]{ctv}, and if $X = \P^d$ then they vanish for $i < d$ \cite[Prop.~3.3(iii)]{ctv}.
\begin{cor} \label{decomp_of_diag}
Suppose that $\dim(X) = d \ge 2$.
\begin{enumerate}
\item If $X$ is stably rational, or more generally if $X$ admits a Chow decomposition of the diagonal \cite[Def.~4.1(i)]{voisin_lueroth}, then $\dimQL(X) \le d - 2$.
\item If $CH_0(X) = \Z$ then there is an integer $N$, independent of $m$, such that $(K/\tau)_d(X,\Z/m)$ is $N$-torsion, and $(K/\tau)_{d-1}(X,\Z/m)$ is $N^2$-torsion; in particular, both vanish when $m$ is relatively prime to $N$.  If $X$ is unirational, we can take $N$ to be the degree of the unirational parametrization.
\end{enumerate}
\end{cor}
\begin{proof}
If $CH_0(X) = \Z$ then $X$ admits a rational Chow decomposition of the diagonal: that is, there is an integer $N$ such that
\begin{equation} \label{rational_decomp}
N\Delta_X = N(X \times x) + Z \text{ in }CH^d(X\times X)
\end{equation}
where $Z$ is a cycle supported on $D \times X$ for some closed subset $D \subsetneq X$.  If $X$ is unirational, we can take $N$ to be the degree of the unirational parametrization.  If $X$ is stably rational, we can take $N=1$; this is called a Chow decomposition of the diagonal.

By \cite[Prop.~3.3]{ctv} and its proof, the existence of a decomposition \eqref{rational_decomp} implies that $H^0(\H^{>0}_X(A))$ and $H^{<d}(\H^d_X(A))$ are $N$-torsion for any coefficient group $A$.  The conclusion follows from taking $A = \Z/m$ and looking at \eqref{K/tau_d} and \eqref{K/tau_d-1}.
\end{proof}
\noindent We remark that the proof of Corollary \ref{K/tau_birational}(a) is not easily adapted to prove Corollary \ref{decomp_of_diag}(a). \medskip

It may be that all rationally connected varieties satisfy $H^0(\H^d_X(\Z/m)) = 0$ for all $m$, or even $H^i(\H^d_X(\Z)) = 0$ for all $i < d$.  For a discussion along these lines, see \cite[\S1.3.2]{voisin_birational}, after \cite[Thm.~1.4]{schreieder_torsion_orders}, and \cite[\S3]{tian_zero_cycles}.\footnote{The ``Kato homology'' $KH_a(X,A)$ appearing in \cite{tian_zero_cycles} is isomorphic to $H^{d-a}(\H^d_X(A))$,
so the conjecture about $H^{<d}(\H^d_X)$ is Tian's statement $\mathbf{RC}(n,k)$ for all $n$ and $k$.}
If this were true, it would simplify and strengthen our Theorem \ref{4fold_thm}, and allow us to say something about rationally connected 5-folds in general. \medskip

We also make the following observation, which is immediate from \eqref{K/tau_d} but seems hard to prove more directly:
\begin{cor} \label{top_unram}
If $X$ and $Y$ are smooth complex projective varieties of dimension $d$, then an equivalence $\Dbcoh(X) \cong \Dbcoh(Y)$ induces an isomorphism of top unramified cohomology groups $H^0(\H^d_X(\Z/m)) \cong H^0(\H^d_Y(\Z/m))$. \qed
\end{cor}
\noindent Since the $m$-torsion part of the Brauer group is isomorphic to $H^0(\H^2_X(\Z/m))$, this sheds some light on the first author's observation \cite{br_not_inv} that the Brauer group of a Calabi--Yau 3-fold is not a derived invariant, in contrast to the Brauer group of a K3 surface which has long been known to be a derived invariant.  In fact Corollary \ref{top_unram} implies that the Brauer group of \emph{any} surface is a derived invariant, which does not seem to have been noticed before. \medskip

Lastly, we repackage a result of Farb, Kisin, and Wolfson and use it to relate Lichtenbaum--Quillen dimension to the outermost Hodge number.\footnote{We thank James Hotchkiss for pointing this out.}

\begin{thm}[Farb, Kisin, and Wolfson \cite{farb-kisin-wolfson}] \label{thm:fkw}
Let $X$ be a smooth complex projective variety.  For all $i \ge 0$ and all primes $p \gg 0$, the rank of the map
\begin{equation} \label{fkw_unramified}
H^i(X,\Z/p) \to H^0(\H^i_X(\Z/p))
\end{equation}
is at least $h^{0,i}(X) = h^i(\O_X)$.
\end{thm}
\begin{proof} 
We apply \cite[Cor.\ 2.2.15]{farb-kisin-wolfson} with $Y = X$ and $D = \varnothing$ to say that for $p \gg 0$, the rank of the map
\begin{equation} \label{fkw_gersten}
H^i_\text{\'et}(X, \Z/p) \to H^i_\text{\'et}(\eta, \Z/p) := \varinjlim_{W \subset X} H^i_{\text{\'et}}(W;\Z/p)
\end{equation}
is at least $h^{i,0}(X) = h^0(\Omega_X^i)$.\footnote{\cite{farb-kisin-wolfson} uses $h^{0,i}$ to denote what everyone else calls $h^{i,0}$, but the two are of course equal.}  By the Artin comparison theorem, we can replace \'etale cohomology with finite coefficients with singular cohomology.  Using the Gersten resolution of $\H^i_X(\Z/p)$ given in \cite{bloch-ogus}, we see that $H^0(\H^i_X(\Z/p))$ injects into $H^i(\eta, \Z/p)$, and the map \eqref{fkw_gersten} factors through it via the map \eqref{fkw_unramified}, so the rank of the latter map is at least $h^{i,0}$ as well.
\end{proof}

\begin{cor}
Let $X$ be a smooth complex projective variety of dimension $d$.  If $h^{0,d}(X) \ne 0$ then $\dimQL(X) = d$.
\end{cor}
\begin{proof}
This is immediate from Theorem \ref{thm:fkw} and \eqref{K/tau_d}.
\end{proof}
We wonder whether the width of the Hodge diamond gives a lower bound on the Lichtenbaum--Quillen dimension in general.  Such a bound could not be sharp, as we will see with Enriques surfaces in \S\ref{enriques}.

\subsection{Construction of the spectral sequence}

To match better with the motivic homotopy theory literature, we switch in this section only to the following setting: fix $m > 1$; let $m' = 2m$ if $m \equiv 2 \pmod 4$, and $m' = m$ otherwise; let $F$ be a field containing $1/m$ and a primitive $m'$$^\text{th}$ root of unity, and with no $m'$-torsion in its Brauer group; and let $X$ be a smooth scheme over $F$.  Then $K_2(F, \Z/m')$ is identified with the roots of unity $\mu_{m'}(F)$, and fixing a primitive root of unity $\zeta$ we can define the spectrum-level Bott maps
\[ \tau\colon  K(X,\Z/m)[2] \to K(X,\Z/m) \]
just as in \S\ref{bott_section}, and consider the cofiber $K/\tau(X,\Z/m)$.  Let
\[ \pi\colon X_\text{\'et} \to X_\text{Zar} \]
be the map from the \'etale site of $X$ to the Zariski site, and let
\[ \H^i_X(\mu_m^{\otimes j}) := R^i \pi_* \mu_m^{\otimes j}. \]
{ \renewcommand{\thethm}{\ref*{spec_seq}$'$}
\addtocounter{thm}{-1}
\begin{thm} \label{spec_seq_etale}
In the setting just described, there is a spectral sequence
\[ E_2^{p,q} = H^p(\H^{-q}_X(\mu_m^{\otimes -q})) \Longrightarrow (K/\tau)_{-p-q}(X,\Z/m), \]
functorial in $X$.
\end{thm}
}
\noindent This implies Theorem \ref{spec_seq} because of the comparison theorem between \'etale cohomology and singular cohomology with finite coefficients.

\begin{proof}
For any Abelian group $A$, we have the motivic cohomology group $H^i_{\text{mot}}(X, A(j))$, which is the $i^\text{th}$ hypercohomology of a certain complex $A(j)$ of sheaves on $X$ in the Zariski topology \cite[Def.~3.4]{mvw}.\footnote{Some readers will be reassured to know that for a smooth separated scheme of finite type over a perfect field, motivic cohomology with $\Z$ coefficients is a re-indexing of Bloch's higher Chow groups: $H^p_\text{mot}(X,\Z(q)) = CH^q(X,2q-p)$.  \cite[Thm.~19.1]{mvw}.}  There is a spectral sequence \cite{fs} \cite{voevodsky_new_approach}
\begin{equation} \label{mot_spec_seq}
E_2^{p,q} = H^{p-q}_\text{mot}(X, A(-q)) \Longrightarrow K_{-p-q}(X, A),
\end{equation}
analogous to the Atiyah--Hirzebruch spectral sequence which computes topological K-theory from singular cohomology.  For $A = \Z/m$, the Beilinson--Lichtenbaum conjecture, which is equivalent to the Bloch--Kato conjecture or norm residue isomorphism theorem, states that the complex $\Z/m(j)$ is quasi-isomorphic to the good truncation\footnote{We are sorry to use $\tau$ for both the Bott maps and the truncation, but both are standard.} $\tau^{\le j} R \pi_*\,\mu_m^{\otimes j}$ \cite[\S1.4]{hw}. \bigskip

Step 1: We argue that our spectrum-level Bott maps \eqref{bott_on_spectra} are compatible with the filtration of $K(X,\Z/m)$ that gives rise to the motivic spectral sequence \eqref{mot_spec_seq}.

We start with Voevodsky's slice filtration
\[ \dotsb \to \Fil^j_\text{mot} K(X) \to \dotsb \to \Fil^1_\text{mot} K(X) \to \Fil^0_\text{mot} K(X) = K(X), \]
proposed in \cite{voevodsky_new_approach} based on some conjectures which have since been settled by Levine in \cite{levine}, with shorter alternate proofs in \cite{bachmann-elmanto}.\footnote{We cannot use Friedlander and Suslin's earlier construction of the motivic spectral sequence \cite{fs} because we need certain multiplicative and functorial properties that it does not provide. We refer to \cite[\S3.1]{e-morrow} for a review.}  The associated graded pieces are the complexes of Abelian groups that compute motivic cohomology:
\[ \text{Cone}(\Fil^{j+1}_\text{mot} K(X) \to \Fil^j_\text{mot} K(X)) =: \gr^j_{\text{mot}} K(X) = \RGamma(X, \Z(j)[2j]). \]
The filtration on $K(X,\Z/m)$ simply comes from tensoring this one with the Moore spectrum $\mathbb S/m$.

When $m \equiv 2 \pmod 4$, we constructed our spectrum-level Bott maps $\tau$ as a composition
\begin{equation} \label{composition}
S^2 \otimes K(X,\Z/m) \xrightarrow{\zeta \otimes 1} K(F,\Z/2m) \otimes K(X,\Z/m) \to K(X,\Z/m),
\end{equation}
where the map $S^2 \to K(F,\Z/2m)$ corresponds to a primitive root of unity $\zeta \in K_2(F,\Z/2m) = \mu_{2m}(F)$, and the second map comes from tensoring the action $K(F) \otimes K(X) \to K(X)$ with the action $\mathbb S/2m \otimes \mathbb S/m \to \mathbb S/m$.  The case $m \not\equiv 2 \pmod 4$ is similar but with $2m$ changed to $m$.

In the first map of \eqref{composition}, we claim that $S^2 \to K(F,\Z/2m)$ lifts uniquely to $\Fil^1_\text{mot} K(F,\Z/2m)$.  To see this, observe that the map
\[ \Fil^1_\text{mot} K(F,\Z/2m) \to \Fil^0_\text{mot} K(F,\Z/2m) \]
induces an isomorphism on $\pi_2$, because its cofiber $\gr^0_\text{mot} K(F,\Z/2m)$ has homotopy groups $\pi_i = H^{-i}_\text{mot}(F,\Z/2m(0))$, and in particular $\pi_2 = \pi_3 = 0$ because there are no negative motivic cohomology groups in this case.

For the second map of \eqref{composition}, it follows from \cite[Thm.~3.6.9]{pelaez} or \cite[\S13.4]{bachmann-hoyois} that the slice filtration is compatible with the ring spectrum structure
\[ K(X) \otimes K(X) \to K(X), \]
so we have
\[ \Fil^i_\text{mot} K(X) \otimes \Fil^j_\text{mot} K(X) \to \Fil^{i+j}_\text{mot} K(X), \]
and in particular
\[ \Fil^1_\text{mot} K(F) \otimes \Fil^{j-1}_\text{mot} K(X) \to \Fil^j_\text{mot} K(X) \]
using the pullback $K(F) \to K(X)$.

Thus our Bott maps are compatible with the filtration:
\[ \tau\colon  \Fil^{j-1}_\text{mot} K(X,\Z/m)[2] \to \Fil^j_\text{mot} K(X,\Z/m). \]

Step 2: Because $\tau$ respects the slice filtration, the cofiber $K/\tau(X,\Z/m)$ inherits a filtration, which gives rise to a spectral sequence computing \linebreak 
$(K/\tau)_*(X,\Z/m)$.  To identify the $E_2$ page of the spectral sequence, we must identify the homotopy groups of the associated graded pieces of the filtration.  These graded pieces fit into exact triangles
\[ \gr^{j-1}_\text{mot} K(X,\Z/m)[2] \to \gr^j_\text{mot} K(X,\Z/m) \to \gr^j_\text{mot} (K/\tau)(X,\Z/m), \]
that is,
\begin{equation} \label{gr_triangle}
\RGamma(X, \Z/m(j-1)[2j]) \to \RGamma(X,\Z/m(j)[2j]) \to \gr^j_\text{mot} (K/\tau)(X,\Z/m).
\end{equation}
We will argue that the first map of \eqref{gr_triangle} is the obvious one, obtained by applying $\RGamma(\underline\quad[2j])$ to the map
\[ \tau^{\le j-1} R \pi_*\,\mu_m^{\otimes{j-1}} \to \tau^{\le j} R \pi_*\,\mu_m^{\otimes j} \]
which includes the shorter truncation into the longer one and multiplies by the root of unity $\zeta \in \mu_m(F)$ chosen earlier, or by $\zeta^2$ if $m \equiv 2 \pmod 4$.  The cone of the latter map is $\H^j_X(\mu_m^{\otimes j})[-j]$, so this will prove that the third term of \eqref{gr_triangle} is $\RGamma(\H^j_X(\mu_m^{\otimes j})[j])$, which gives the $E_2$ page promised in the theorem.

The multiplication map $K(X) \otimes K(X) \to K(X)$ is compatible with the slice filtration, as we have said, and it follows from \cite[Thm.~3.6.23]{pelaez} that the resulting maps on associated graded pieces coincide (after a shift) with the usual multiplication
\[ \RGamma(X, \Z(i)) \otimes \RGamma(X, \Z(j)) \to \RGamma(X, \Z(i+j)) \]
in motivic cohomology, constructed for example in \cite[Constr.~3.11]{mvw}.  Tensoring with the action $\mathbb S/2m \otimes \mathbb S/m \to \mathbb S/m$, we find that the maps on associated graded pieces induced by the action
\[ K(F,\Z/2m) \otimes K(X,\Z/m) \to K(X,\Z/m) \]
coincide (again after a shift) with the motivic multiplication
\[ \RGamma(F,\Z/2m(i)) \otimes \RGamma(X,\Z/m(j)) \to \RGamma(X,\Z/m(i+j)). \]
Reviewing our construction of the spectrum-level Bott maps, we find that the first map of \eqref{gr_triangle} is given by motivic multiplication with
\[ \zeta \in H^0_\text{mot}(F,\Z/2m(1)) = \mu_{2m}(F). \]
By \cite[\S5]{weibel_products}, 
the multiplication in motivic cohomology is compatible, via the maps
\begin{equation} \label{motivic_to_etale}
\Z/m(i) \to R\pi_* \mu_m^{\otimes i},
\end{equation}
with the obvious multiplication
\begin{equation} \label{etale_mult}
R\pi_* \mu_m^{\otimes i} \otimes R\pi_* \mu_m^{\otimes j} \to R\pi_* \mu_m^{\otimes i+j}.
\end{equation}
Moreover, the maps \eqref{motivic_to_etale} identify $\Z/m(i)$ with $\tau^{\le i} R\pi_* \mu_m^{\otimes i}$ by \cite[Thm.~C]{hw}, and any map of the form \eqref{etale_mult} must respect truncations and take \linebreak 
$\tau^{\le i} R\pi_* \mu_m^{\otimes i} \otimes \tau^{\le j} R\pi_* \mu_m^{\otimes j}$ into $\tau^{\le i+j} R\pi_* \mu_m^{\otimes i+j}$.  Everything is similar with the first $m$ changed to $2m$.  Thus the first map of \eqref{gr_triangle} is what we claimed.
\end{proof}

\subsection{Conic bundles and twisted K-theory} \label{conic_section}

As a first application of our motivic spectral sequence, we revisit Corollary \ref{Pn_is_QL0}, which said that if $P \to X$ is a $\P^r$ bundle then $\dimQL(P) \ge \dimQL(X)$, with equality if the bundle is locally trivial in the Zariski topology.  Whether equality holds in general is a surprisingly stubborn question.  Here we prove only that if $r=1$ then $\dimQL(P) \le \dim(X)$.
\begin{prop} \label{conic_fibr}
If a smooth complex projective variety $Y$ admits a map $f\colon Y \to Z$ whose generic fiber is a conic, then $\dimQL(Y) \le \dim(Z)$.
\end{prop}
\begin{proof}
Let $d = \dim(Y)$.  We claim that $H^0(\H^d_Y(\Q/\Z)) = 0$: when $d=3$, Colliot-Th\'el\`ene and Voisin \cite[Cor.~8.2]{ctv} proved this using Kahn, Rost, and Sujatha's \cite[Thm.~5]{krs1}, and when $d$ is arbitrary, it follows in the same way from \cite[Prop.~A.1]{krs1}.\footnote{See \cite[Lem.~9]{fu_tian} for another version of the same argument.}  
Next, we have exact sequences
\[ 0 \to \H^j_Y(\Z/m) \to \H^j_Y(\Q/\Z) \xrightarrow{\cdot m} \H^j_Y(\Q/\Z) \to 0 \]
for all $j$ and all $m$; this follows for example from \cite[Cor.~1.5]{voisin_deg_4}.  Thus we find that $H^0(\H^d_Y(\Z/m)) = 0$, so \eqref{K/tau_d} gives $(K/\tau)_d(Y,\Z/m) = 0$, so $\dimQL(Y) \le d-1 = \dim(Z)$.
\end{proof}

Hotchkiss \cite[Thm.~6.1]{hotchkiss} proved that if $X$ is an Abelian 3-fold and $P \to X$ is an \'etale $\P^3$ bundle whose Brauer class has order 2 but cannot be represented by an \'etale $\P^1$ bundle, then the integral Hodge conjecture fails for $H^6(P,\Z)$. He constructed a similar example over products of curves in \cite[Thm.~6.4]{hotchkiss}.  These may be good places to look for counterexamples to our expectation that $\dimQL(P) = \dimQL(X)$.\medskip

For a Brauer class $\alpha \in \Br(X)$, let $\dimQL(X,\alpha)$ denote the Lichtenbaum--Quillen dimension of the derived category of twisted sheaves $\Dbcoh(X,\alpha)$.  The question above amounts to asking whether $\dimQL(X,\alpha) \le \dimQL(X)$ for all $\alpha$.  To rephrase what we have already said in this language:
\begin{cor} \label{twisted}
Let $X$ be a smooth complex projective variety, and let $\alpha \in \Br(X)$ be a Brauer class of index 2.\footnote{Meaning that $\alpha$ can be represented by an \'etale $\P^1$ bundle, or equivalently by an Azumaya algebra of rank 4.  De Jong \cite{de_jong} proved that if $\dim(X) = 2$ then every 2-torsion class has index 2.}  Then $\dimQL(X,\alpha) \le \dim(X)$. \qed
\end{cor}
\noindent But funny things can happen: in \S\ref{enriques} and \S\ref{AM_3folds} we will see examples where $\dim(X) = 2$ or 3 and $\dimQL(X) = 2$, but $\dimQL(X,\alpha) = 0$.

\subsection{Two results of Pedrini and Weibel}

Pedrini and Weibel's purpose in \cite{pw_surfaces} and \cite{pw_varieties} was to use the Lichtenbaum--Quillen conjecture to describe the K-theory of smooth complex varieties, not just with finite coefficients but with integral coefficients.  The following proposition, which we will use in \S\ref{fano_3folds}, \S\ref{cubic_gm_sec}, and \S\ref{cubic_5fold_sec}, is a reworking of their \cite[Example 6.7]{pw_surfaces} and \cite[Thm.~C]{pw_varieties}, and goes back to a remark of Thomason \cite[after eq.~(0.5)]{thomason_etale}.  The point is that the Lichtenbaum--Quillen conjecture only sheds light on $K_n(X)$ for $n \ge \dim(X) - 2$, leaving small $n$ in question, but the present paper shows that many higher-dimensional varieties behave like lower-dimensional ones from the perspective of the Lichtenbaum--Quillen conjecture.

\begin{prop} \label{pw_prop}
If $X$ is a smooth complex projective variety, then
\[ K_n(X) = \operatorname{tors}(KU^{-n}(X)) \oplus (\Q/\Z)^{\rank KU^{-n-1}(X)} \oplus K_n(X,\Q) \]
for all $n \ge \max(1,\, \dimQL(X)-2)$.  The same holds if $X$ is replaced by an admissible subcategory $\A \subset \Dbcoh(X)$.
\end{prop}
\begin{proof}
For brevity set $d = \dimQL(X)$.  By definition the comparison maps $K_n(X,\Z/m) \to KU^{-n}(X,\Z/m)$ are isomorphisms for $n \ge d-1$ and injective for $n = d-2$.  Taking direct limits over $m$, we find that the comparison maps $K_n(X,\Q/\Z) \to KU^{-n}(X,\Q/\Z)$ are isomorphisms for $n \ge d-1$ and injective for $n = d-2$.\footnote{The reader may wonder why we did not define $\dimQL(X)$ in terms of the comparison maps with $\Q/\Z$ coefficients in the first place, rather than $\Z/m$ coefficients for all $m$.  The reason is that there are no Bott maps with $\Q/\Z$ coefficients, so we would lose the spectral sequence of Theorem \ref{spec_seq}.}

By \cite[\S6.3]{gillet} the image of the comparison map $K_n(X) \to KU^{-n}(X)$ is torsion for all $n \ge 1$, so the comparison maps $K_n(X,\Q) \to KU^{-n}(X,\Q)$ are zero for all $n \ge 1$.

Thus for all $n \ge \max(1,\,d-2)$ we have a commutative diagram with exact rows
\[ \centerline{ \xymatrix@C=1.5em{
K_{n+1}(X,\Q) \ar[d]^-0 \ar[r] & K_{n+1}(X,\Q/\Z) \ar@{=}[d] \ar[r] & K_n(X) \ar[r] \ar[d] & K_n(X,\Q) \ar[r] \ar[d]^0 & K_n(X,\Q/\Z) \ar@{^(->}[d] \\
KU^{-n-1}(X,\Q) \ar[r] & KU^{-n-1}(X,\Q/\Z) \ar[r] & KU^{-n}(X) \ar[r] & KU^{-n}(X,\Q) \ar[r] & KU^{-n}(X,\Q/\Z).
} } \]
In the top row, we find that the map $K_n(X,\Q) \to K_n(X,\Q/\Z)$ must be zero, and similarly with $K_{n+1}$, so we get a short exact sequence
\begin{equation} \label{to_split}
0 \to K_{n+1}(X,\Q/\Z) \to K_n(X) \to K_n(X,\Q) \to 0.
\end{equation}
From the bottom row we get a short exact sequence
\[ 0 \to (\Q/\Z)^{\rank KU^{-n-1}(X)} \to KU^{-n-1}(X,\Q/\Z) \to \operatorname{tors}(KU^{-n}(X)) \to 0, \]
which must split because $\Q/\Z$ is divisible and hence injective.  Thus \eqref{to_split} splits because $\Ext^1_\Z(\Q,\Z/m) = 0$ and $\Ext^1_\Z(\Q,\Q/\Z) = 0$. 
\end{proof}

The proof above might fail for a general smooth proper dg-category $\mathcal C$, because we do not know whether the image of the comparison map $K_n(\mathcal C) \to KU^{-n}(\mathcal C)$ is torsion for $n \ge 1$.  But for completeness we record the following proposition, which is a reworking of \cite[Cor.~4.3 and 6.4]{pw_varieties}.

\begin{prop}
If $X$ is a smooth complex projective variety, then the image of the comparison map
\[ \eta_n\colon K_n(X) \to KU^{-n}(X) \]
contains the torsion subgroup of $KU^{-n}(X)$ for $n \ge \dimQL(X)-2$.
The same holds if $X$ is replaced by an admissible subcategory $\A \subset \Dbcoh(X)$, or an arbitrary $\C$-linear dg-category $\mathcal C$.
\end{prop}
\begin{proof}
Take the spectrum-level comparison map $\eta\colon K(X) \to KU(X)$, and let $\cone \eta$ be its cofiber.  It is enough to prove that $\pi_n(\cone \eta)$ is torsion-free for $n \ge d - 2$, where $d = \dimQL(X)$.

Consider the exact triangle
\[ \cone \eta \xrightarrow {\cdot m} \cone \eta \to (\cone \eta)/m. \]
We want to prove that the first map induces an injection on $\pi_n$ for all $m$ and $n \ge d - 2$, so it is enough to prove that the third term has $\pi_n = 0$ for $n \ge d - 1$.

By the $3 \times 3$ lemma we get an exact triangle
\[ K(X,\Z/m) \xrightarrow\eta KU(X,\Z/m) \to (\cone \eta)/m. \]
By hypothesis, the first map induces an isomorphism on $\pi_n$ for $n \ge d - 1$ and an injection for $n = d - 2$, so the third term has $\pi_n = 0$ for $n \ge d - 1$, as desired.
\end{proof}

%% file: surfaces.tex

\section{Surfaces} \label{surfaces}

Castelnuovo proved that a smooth complex projective surface $X$ is rational if and only if the irregularity $q := h^1(\O_X)$ and the second plurigenus $P_2 := h^0(\omega_X^{\otimes 2})$ both vanish.  But it is not enough for $q$ and the geometric genus $p_g = P_1 = h^0(\omega_X)$ to vanish, as Enriques surfaces demonstrate, nor is it enough to have $p_g = q = 0$ and $\pi_1 = 0$, as the surfaces of Dolgachev \cite{dolgachev} and Barlow \cite{barlow} demonstrate.

In this section we study the Lichtenbaum--Quillen dimension of some surfaces, not to shed light on the rationality problem, but to get a feeling for how to use our spectral sequence, and because they will come up later in connection with 3-folds and 4-folds.

\subsection{Analysis of the spectral sequence} \label{spec_seq_surfaces}

The spectral sequence from Theorem \ref{spec_seq} relates $(K/\tau)_*(X,\Z/m)$ to the groups $H^i(\H^j_X(\Z/m))$.  We will get our hands on the latter using two tools.

First, we will relate them to $H^i(\H^j_X(\Z))$ using a result of Colliot-Th\'el\`ene and Voisin \cite[Thm.~3.1]{ctv}, which they deduced from the norm residue isomorphism theorem using an argument of Bloch and Srinivas \cite[Thm.~1(ii)]{bloch-srinivas} \cite[Rmk.~5.8]{bloch_book}: the sheaves $\H^j_X(\Z)$ are torsion-free, so the long exact sequence in $\H^*_X$ coming from $0 \to \Z \xrightarrow{\cdot m} \Z \to \Z/m \to 0$ splits into short exact sequences
\[ 0 \to \H^j_X(\Z) \xrightarrow{\cdot m} \H^j_X(\Z) \to \H^j_X(\Z/m) \to 0, \]
which in turn give long exact sequences
\begin{multline} \label{ctv_les}
0 \to H^0(\H^j_X(\Z)) \xrightarrow{\cdot m} H^0(\H^j_X(\Z)) \to  H^0(\H^j_X(\Z/m)) \\
\to H^1(\H^j_X(\Z)) \xrightarrow{\cdot m} H^1(\H^j_X(\Z)) \to  H^1(\H^j_X(\Z/m)) \to \dotsb
\end{multline}
for each $j$.

Second, we use the Bloch--Ogus--Leray spectral sequence
\begin{equation} \label{bol_ss}
E_2^{p,q} = H^p(\H^q_X(A)) \Longrightarrow H^{p+q}_\text{sing}(X,A).
\end{equation}
The resulting filtration on $H^*_\text{sing}(X,A)$ is the coniveau filtration.  The $E_2$ page looks like \eqref{our_E_2} above, but upside down (so the differentials go a different way),
and we will use $\Z$ coefficients rather than $\Z/m$.  In addition to the vanishing of $H^i(\H^j_X(A))$ for $i > j$ or $j > d$, we will need \cite[Cor.~7.4]{bloch-ogus} which says that $H^i(\H^i_X(\Z))$ is the quotient of $CH^i(X)$ by algebraic equivalence, and the resulting edge map to $H^{2i}(X,\Z)$ is the cycle class map.

{ \renewcommand{\thethm}{\ref{surface_thm}}
\addtocounter{thm}{-1}
\begin{thm}
Let $X$ be a smooth complex projective surface.
\begin{enumerate}
\item $\dimQL(X) \le 1$ if and only if the geometric genus $p_g = 0$, and $H^3(X,\Z)$ is torsion-free.
\item $\dimQL(X) = 0$ if and only if in addition, the irregularity $q = 0$.
\end{enumerate}
\end{thm}
}

\begin{proof}
We may assume that $X$ is connected.  We see that the $E_2$ page of the Bloch--Ogus--Leray spectral sequence \eqref{bol_ss} with $\Z$ coefficients looks like this:\footnote{For the benefit of any topologists reading, we recall that the N\'eron--Severi group $\NS(X)$ is the image of $c_1\colon \Pic(X) \to H^2(X,\Z)$.}
\pagebreak 
\[ \xymatrix@=1ex{
2 & \dfrac{H^2(X,\Z)}{\NS(X)} & H^3(X,\Z) & \Z \\
1 & H^1(X,\Z) & \NS(X) \\
q=0 & \Z \\
\ar@{-}[]+<1.6em,.8em>;[rrr]+<1.6em,.8em> 
\ar@{-}[]+<1.6em,.8em>;[uuu]+<1.6em,.8em> 
& p=0 & 1 & 2
} \]

(a) From \eqref{K/tau_d} we have
\[ (K/\tau)_2(X,\Z/m) = H^0(\H^2_X(\Z/m)). \]
From the long exact sequence \eqref{ctv_les} we see that this fits into a short exact sequence
\[ 0 \to \dfrac{H^2(X,\Z)}{\NS(X)} \otimes \Z/m \to H^0(\H^2_X(\Z/m)) \to m\text{-tors}(H^3(X,\Z)) \to 0. \]
Because $H^2(X,\Z)/\NS(X)$ is finitely generated, the first term vanishes for all $m$ if and only if $H^2(X,\Z)/\NS(X) = 0$, and by the Lefschetz theorem on $(1,1)$ classes, this happens if and only if $p_g = 0$.  Thus we see that $(K/\tau)_2(X,\Z/m) = 0$ for all $m$ if and only if $p_g = 0$ and $H^3(X,\Z)$ is torsion-free. \medskip

(b) Now $\NS(X) = H^2(X,\Z)$ is torsion-free by Poincar\'e duality and the universal coefficient theorem, and by Hodge theory we have $H^1(X,\Z) = \Z^{2q}$ and $H^3(X,\Z) = \Z^{2q}$, so from the long exact sequence \eqref{ctv_les} we find that $H^1(\H^2_X(\Z/m)) = (\Z/m)^{2q}$ and $H^0(\H^1_X(\Z/m)) = (\Z/m)^{2q}$.  Thus the exact sequence \eqref{K/tau_d-1} becomes
\[ 0 \to (\Z/m)^{2q} \to (K/\tau)_1(X,\Z/m) \to (\Z/m)^{2q} \to \Z/m, \]
from which we see that $(K/\tau)_1(X,\Z/m) = 0$ if and only if $q=0$.
\end{proof}

\subsection{Enriques surfaces} \label{enriques}

An Enriques surface $X$ is the quotient of a K3 surface by a free involution.  It has $p_g = q = 0$ but $H_1(X,\Z) = \Z/2$, so $\dimQL(X) = 2$: that is, the map $K_n(X,\Z/m) \to KU^{-n}(X,\Z/m)$ is injective but not surjective for $n=0$ and some $m$, in fact $m=2$.  Let us take a moment to see what the maps look like explicitly.

The integral cohomology of an Enriques surface is
\[ H^i(X,\Z) = \begin{cases} 
\Z & i = 0 \\
0 & i = 1 \\
\Z^{10} \oplus \Z/2 & i = 2 \\
\Z/2 & i = 3 \\
\Z & i = 4,
\end{cases} \]
and the Atiyah--Hirzebruch spectral sequence gives
\begin{align*}
KU^\text{even}(X) &= \Z^{12} \oplus \Z/2 & KU^\text{odd}(X) &= \Z/2.
\end{align*}
Most of the details are given in \cite[\S1]{aw}; to finish the calculation of $KU^0$, one observes that the edge map $H^4(X,\Z) \hookrightarrow KU^0(X)$ is split by the pushforward $KU^0(X) \to KU^0(\text{point})$.

We have a diagram of long exact sequences
\begin{equation} \label{les_KU}
\xymatrix{
K_0(X) \ar@{->>}[d] \ar[r]^2 & K_0(X) \ar@{->>}[d] \ar[r] & K_0(X, \Z/2) \ar@{^(->}[d] \ar[r] & 0 \ar[d] \ar[r] & 0 \ar[d] \\
KU^0(X) \ar[r]^2 & KU^0(X) \ar[r] & KU^0(X, \Z/2) \ar[r] & KU^1(X) \ar[r]^2 & KU^1(X),
}
\end{equation}
where the surjectivity of the vertical arrows on the left comes from Lemma \ref{K-KU} below, and the injectivity of the middle vertical arrow comes from the Lichtenbaum--Quillen conjecture.  This yields a diagram of short exact sequences
\[ \xymatrix{
0 \ar[r] & K_0(X)/2 \ar@{=}[d] \ar[r] & K_0(X, \Z/2) \ar@{^(->}[d] \ar[r] & 0 \ar[d] \ar[r] & 0  \\
0 \ar[r] & KU^0(X)/2 \ar[r] & KU^0(X, \Z/2) \ar[r] & 2\text{-tors}(KU^1(X)) \ar[r] & 0.
} \]
Now $KU^0(X, \Z/2)$ is a $\Z/2$-vector space by \cite[Thm.~2.3]{at1},
so the second row must read
\begin{equation} \label{13-14-1}
0 \to (\Z/2)^{13} \to (\Z/2)^{14} \to \Z/2 \to 0.
\end{equation}
Thus the comparison map $K_0(X,\Z/2) \to KU^0(X,\Z/2)$ is identified with the first map of \eqref{13-14-1}.  This is indeed injective and not surjective, and we have seen how the cokernel came from the torsion in $H^3(X,\Z)$. \bigskip

We conclude with the surprising twisted example promised in \S\ref{conic_section}.  We wonder whether $\Dbcoh(X,\alpha)$ admits a full exceptional collection.

\begin{prop} \label{twisted_enriques}
If $X$ is a (complex) Enriques surface and $\alpha$ is the non-zero element of $\Br(X) = H^3(X,\Z) = \Z/2$, then $\dimQL(X,\alpha) = 0$.
\end{prop}
\begin{proof}
Because $X$ is a surface and $\alpha$ is 2-torsion, $\alpha$ is represented by an \'etale $\P^1$ bundle $\pi\colon P \to X$ by \cite{de_jong}.  We saw in \S\ref{semi-orth} that
\[ K_*(P) = K_*(X) \oplus K_*(X,\alpha), \]
and similarly with $KU$ and with $\Z/m$ coefficients.

Corollary \ref{twisted} gives $\dimQL(X,\alpha) \le 2$, so to prove that $\dimQL(X,\alpha) = 0$ we must show that the comparison map $K_0(X,\alpha,\Z/m) \to KU^0(X,\alpha,\Z/m)$ is surjective as well as injective, and that $KU^1(X,\alpha,\Z/m) = 0$.

We see that the 3-fold $P$ satisfies the hypotheses of Lemma \ref{K-KU} below, so the comparison map $K_0(P)/m \to KU^0(P)/m$ is surjective, so $K_0(X,\alpha)/m \to KU^0(X,\alpha)/m$ is surjective.  From the $\alpha$-twisted analogue of \eqref{les_KU}, we see it is enough to prove that $KU^1(X,\alpha) = 0$ and $KU^2(X,\alpha)$ is torsion-free.

From the Gysin sequence
\[ \dotsb \to H^i(X,\Z) \xrightarrow{\pi^*} H^i(P,\Z) \xrightarrow{\pi_*} H^{i-2}(X,\Z) \xrightarrow{\cdot\alpha} H^{i+1}(X,\Z) \to \dotsb \]
we find that
\[ H^i(P,\Z) = \begin{cases} 
\Z & i = 0 \\
0 & i = 1 \\
\Z^{11} \oplus \Z/2 & i = 2 \\
0 & i = 3 \\
\Z^{11} & i = 4 \\
\Z/2 & i = 5 \\
\Z & i = 6.
\end{cases} \]
Here the only tricky step is in determining that the extension
\[ 0 \to \Z \to H^4(P,\Z) \to \Z^{10} \oplus \Z/2 \to 0 \]
is not split; this is because the torsion in $H^4(P,\Z)$ is isomorphic to the torsion in $H_3(P,\Z) = H^3(P,\Z) = 0$.

The Atiyah--Hirzebruch spectral sequence for $P$ degenerates at the $E_2$ page, and the filtration splits: in running the spectral sequence and solving the extension problem, we can only lose torsion, but $KU^*(P)$ contains $KU^*(X)$ as a direct summand, so all the torsion must survive by our earlier calculation of $KU^*(X)$.  Thus
\begin{align*}
KU^\text{even}(P) &= \Z^{24} \oplus \Z/2 & KU^\text{odd}(P) &= \Z/2,
\end{align*}
and thus
\begin{align*}
KU^\text{even}(X,\alpha) &= \Z^{12} & KU^\text{odd}(X,\alpha) &= 0
\end{align*}
as desired.  This agrees with the calculation of $KU^\text{odd}(X,\alpha)$ in \cite[\S1]{aw} using the twisted Atiyah--Hirzebruch spectral sequence.
\end{proof}

\begin{lem} \label{K-KU}
Let $X$ be a smooth complex projective surface with $h^{2,0} = 0$, or a smooth uniruled 3-fold with $h^{2,0} = 0$.  Then the comparison map $K_0(X) \to KU^0(X)$ is surjective.
\end{lem}
\begin{proof}
The cycle class map $CH^i(X) \to H^{2i}(X,\Z)$ is surjective for all $i$: for $i=0$ and $i = \dim(X)$ this is clear, for $i=1$ it is the Lefschetz $(1,1)$-theorem, and for uniruled 3-folds and $i=2$, it is Voisin's \cite[Thm.~2]{voisin_uniruled}.

Now take the filtration by codimension of support,
\[ \dotsb \subset F^2 \subset F^1 \subset F^0 = K_0(X). \]
By \cite[Ex.~15.1.5 and 15.3.6]{fulton}, there is a surjection
\[ CH^i(X) \twoheadrightarrow F^i/F^{i+1} \]
for each $i$.  When $X$ is a surface, the Atiyah--Hirzebruch spectral sequence collapses by \cite[proof of Prop.~1.1]{aw}, and when $X$ is a 3-fold it collapses by a similar argument, so it gives a filtration on $KU^0(X)$ whose associated graded pieces are $H^{2i}(X,\Z)$.  The comparison map $K_0(X) \to KU^0(X)$ is compatible with the filtrations, and the resulting maps $CH^i(X) \to H^{2i}(X,\Z)$ are the cycle class maps \cite[Thm.~1.3]{fhw}.  Since these are surjective, we find that $K_0(X) \to KU^0(X)$ is surjective.
\end{proof}

\subsection{Barlow and Dolgachev surfaces}

Barlow surfaces were the first known surfaces of general type with $p_g = q = 0$ and $H_1(X,\Z) = 0$.\footnote{More recent examples are collected in \cite[Table 1]{gen_type_pg_zero}.  Their derived categories do not seem to have been studied.}  Thus they satisfy $\dimQL(X) = 0$ but are not rational, so Lichtenbaum--Quillen dimension does not give a complete obstruction to rationality even in dimension 2.

B\"ohning, von Bothmer, Katzarkov, and Sosna \cite{bvbks} showed that the derived category of a generic determinantal Barlow surface admits a semi-orthogonal decomposition
\[ \Dbcoh(X) = \langle \A_X, L_1, \dotsc, L_{11} \rangle, \]
where $L_1, \dotsc, L_{11}$ are exceptional line bundles and $\A_X$ is a ``phantom category,'' meaning that it is non-zero but its Hochschild homology and $K_0$ vanish.  In fact the K-motive of $\A_X$ vanishes, as Sosna showed in \cite[Cor.~4.9]{sosna_remarks}, so $K_n(\A_X) = 0$ for all $n$.  But the Hochschild cohomology of $\A_X$ does not vanish.

We can re-prove that $K_n(\A_X) = 0$ for all $n$ as follows; the hypothesis $\dimQL(\A_X) \le 0$ follows from Proposition \ref{sod_vs_QL}.
\begin{prop}
Let $X$ be a smooth complex projective variety and $\A \subset \Dbcoh(X)$ be an admissible subcategory.  If $K_0(\A) = 0$ and $\dimQL(\A) \le 0$, then $K_n(\A) = 0$ for all $n$.
\end{prop}
\begin{proof}
The hypothesis $\dimQL(\A) \le 0$ gives
\begin{align*}
K_0(\A,\Z/m) = K_2(\A,\Z/m) = K_4(\A,\Z/m) = \dotsb \\
K_{-1}(\A,\Z/m) = K_1(\A,\Z/m) = K_3(\A,\Z/m) = \dotsb
\end{align*}
for all $m$.  But $K_n(\A,\Z/m) = 0$ for $n < 0$, and $K_0(\A,\Z/m) = K_0(\A)/m = 0$ by hypothesis, so $K_n(\A,\Z/m) = 0$ for all $n$.  From the long exact sequence
\[ \dotsb \to K_{n+1}(\A,\Z/m) \to K_n(\A) \xrightarrow{\cdot m} K_n(\A) \to K_n(\A,\Z/m) \to \dotsb \]
we see that $K_n(\A)$ is a torsion-free divisible group for all $n$, which implies that $K_n(\A) = K_n(\A,\Q)$.  On the other hand, Gorchinskiy and Orlov proved in \cite[Thm.~5.5]{go} that $K_0(\A,\Q) = 0$ implies $K_n(\A, \Q) = 0$ for all $n$.
\end{proof}

The same analysis applies to Dolgachev surfaces, which are surfaces of Kodaira dimension 1 with $p_g = q = 0$ and $H_1(X,\Z) = 0$.  Cho and Lee proved in \cite[Thm.~1.4]{cho_lee} that for well-chosen parameters, the derived category decomposes into an exceptional collection of 12 line bundles and a phantom category.

%% file: 3folds.tex

\section{Rationally connected 3-folds} \label{3-folds}

Because we were originally motivated by questions of rationality, we restrict our attention to rationally connected varieties from here on.  This simplifies our analysis by killing some unramified cohomology groups.

\subsection{Analysis of the spectral sequence}

The following result uses the same ingredients as \cite[Prop.~6.3]{ctv}, but we give all the details as a model for our study of 4-folds in \S\ref{4fold_thm_sec}.

{ \renewcommand{\thethm}{\ref{3fold_thm}}
\addtocounter{thm}{-1}
\begin{thm}
Let $X$ be a smooth complex projective 3-fold that is rationally connected.
\begin{enumerate}
\item $\dimQL(X) \le 2$.
\item $\dimQL(X) \le 1$ if and only if $H^3(X,\Z)$ is torsion-free.
\item $\dimQL(X) = 0$ if and only if $H^3(X,\Z) = 0$.
\end{enumerate}
\end{thm}
}

\begin{proof}
First we argue that the $E_2$ page of the Bloch--Ogus--Leray spectral sequence \eqref{bol_ss} with $\Z$ coefficients looks like this:
\[ \xymatrix@=1ex{
3 & 0 & 0 & 0 & \Z \\
2 & 0 & H^3(X,\Z) & \dfrac{CH^2(X)}{\text{\footnotesize alg.\,equiv.}} \\
1 & 0 & \NS(X) \\
q=0 & \Z \\
\ar@{-}[]+<1.6em,.8em>;[rrrr]+<1.6em,.8em> 
\ar@{-}[]+<1.6em,.8em>;[uuuu]+<1.6em,.8em> 
& p=0 & 1 & 2 & 3
} \]

In the column $p=0$, we have $H^0(\H^q_X(\Z)) = 0$ for $q > 0$ by \cite[Prop.~3.3(i)]{ctv}, because $CH_0(X) = \Z$.  Thus the differentials vanish.

The diagonal $p+q = 4$ gives a 2-step filtration of $H^4(X,\Z)$ where the sub-object is $H^2(\H^2_X(\Z))$, which we have seen is $CH^2(X)$ modulo algebraic equivalence, so the quotient $H^1(\H^3_X(\Z))$ is the cokernel of the cycle class map $CH^2(X) \to H^4(X,\Z)$.  But a rationally connected variety is uniruled \cite[Cor.~4.11 and Cor.~4.17(b)]{debarre_book}, so this cycle class map is surjective by Voisin's \cite[Thm.~2]{voisin_uniruled}.

On the diagonal $p+q=5$, we have $H^5(X,\Z) = 0$ because rationally connected varieties are simply connected \cite[Cor.~4.18(c)]{debarre_book}.

Thus we have justified the $E_2$ page given above. \medskip

(a) From \eqref{K/tau_d} we have
\[ (K/\tau)_3(X,\Z/m) = H^0(\H^3_X(\Z/m)), \]
and from the long exact sequence \eqref{ctv_les} we see that this vanishes. \medskip

(b) The exact sequence \eqref{K/tau_d-1} reads
\begin{multline*}
0 \to H^1(\H^3_X(\Z/m)) \to (K/\tau)_2(X, \Z/m) \\
\to H^0(\H^2_X(\Z/m)) \to H^2(\H^3_X(\Z/m)). \end{multline*}
From the long exact sequence \eqref{ctv_les} and the $E_2$ page above we find that $H^1(\H^3_X(\Z/m)) = 0$ and $H^2(\H^3_X(\Z/m)) = 0$, so
\[ (K/\tau)_2(X,\Z/m) = H^0(\H^2_X(\Z/m)) = m\text{-tors}(H^3(X,\Z)). \]
Thus $(K/\tau)_2(X,\Z/m) = 0$ for all $m$ if and only if $H^3(X,\Z)$ is torsion-free. \medskip
 
(c) Now the comparison maps $K_n(X,\Z/m) \to KU^{-n}(X,\Z/m)$ are isomorphisms for $n \ge 0$ and injective for $n = -1$.  We always have $K_{<0} = 0$, so we will have $\dimQL(X) = 0$ if and only if $KU^1(X,\Z/m) = 0$ for all $m$.  We know that $H_1(X,\Z) = 0$ and $H^3(X,\Z)$ is torsion-free, so all of $H^*(X,\Z)$ is torsion-free by Poincar\'e duality and the universal coefficient theorem.  Using the Atiyah--Hirzebruch spectral sequence we find that $KU^*(X)$ is torsion-free and $KU^1(X) = H^3(X,\Z)$, so $KU^1(X,\Z/m) = KU^1(X)/m$ vanishes if and only if $H^3(X,\Z) = 0$.
\end{proof}

\subsection{Artin--Mumford 3-folds} \label{AM_3folds}

Artin and Mumford \cite{am} constructed a unirational 3-fold $X$ with $H^3(X,\Z) = \Z/2$ by taking the quartic surface cut out by the determinant of a general $4 \times 4$ symmetric matrix of linear forms, taking the double cover of $\P^3$ branched over that surface, and blowing up the 10 ordinary double points.  By Theorem \ref{3fold_thm} we have $\dimQL(X) = 2$.

Again this means that the map $K_n(X,\Z/2) \to KU^{-n}(X,\Z/2)$ is injective but not surjective for $n=0$, and it is an interesting exercise to see so directly as we did for Enriques surfaces in \S\ref{enriques}.  But the similarity between Artin--Mumford 3-folds and Enriques surfaces is no coincidence: Li, Nuer, Stellari, and Zhao \cite[Thm.~B]{lnsz}, building on partial results by Ingalls and Kuznetsov \cite[Thm.~4.3]{ik} and Hosono and Takagi \cite[\S1.2]{ht_am}, proved that there is a semi-orthogonal decomposition
\[ \Dbcoh(X) = \langle \Dbcoh(S), E_1, \dotsc, E_{12} \rangle \]
where $S$ is an Enriques surface and $E_1, \dotsc, E_{12}$ are exceptional sheaves.  Thus the fact that $\dimQL(X) = \dimQL(S) = 2$ agrees with Proposition \ref{sod_vs_QL}.

The following result echoes Proposition \ref{twisted_enriques}.  We wonder whether the twisted derived categories of Artin--Mumford 3-folds and Enriques surfaces are related, as the untwisted derived categories are.
\begin{prop}
If $X$ is an Artin--Mumford 3-fold and $\alpha$ is the non-zero element of $\Br(X) = H^3(X,\Z) = \Z/2$, then $\dimQL(X,\alpha) = 0$.  
\end{prop}
\begin{proof}
We will argue that $\Dbcoh(X,\alpha)$ admits a full exceptional collection.

Let $\bar X$ be the double cover of $\P^3$ branched over the determinantal quartic surface mentioned above, so $\bar X$ has 10 ordinary double points.  Let $\hat X$ be a small resolution of $\bar X$, obtained by replacing each ordinary double point by a $\P^1$ in one of two ways; this is only a Moishezon manifold or algebraic space, not a projective variety \cite[Prop.~4.0.4]{nick_thesis}, but this will not cause any trouble.  The big resolution $X$ is the blow-up of $\hat X$ along the 10 exceptional $\P^1$s.  Alternatively, we could have constructed $\hat X$ from $X$ by contracting each exceptional $\P^1 \times \P^1 \subset X$ down to a $\P^1$ in one of two ways.

A Brauer class on the small resolution $\hat X$ was constructed in the first author's thesis \cite[\S4.1.3]{nick_thesis}; here we will call it $\hat\alpha \in \Br(\hat X)$.  Cheng and Olander \cite[Thm.~A]{cheng-olander}, refining and extending the results of \cite{nick_thesis}, proved that $\Dbcoh(\hat X, \hat\alpha)$ is equivalent to Kuznetsov's category $\Dbcoh(\P^3, \mathcal B_0)$, where $\mathcal B_0$ is a sheaf of even parts of Clifford algebras.  Kuznetsov \cite[Thm.~5.5]{kuznetsov_quadrics} proved that $\Dbcoh(\P^3, \mathcal B_0)$ has a full exceptional collection of length 4; to extract this statement from the reference, note that the intersection of four quadrics in $\P^3$ is empty.

Let $\pi$ be the map $X \to \hat X$, which blows up 10 copies of $\P^1$.  A straightforward extension of Example \ref{sod_examples}(d) gives a semi-orthogonal decomposition
\[ \Dbcoh(X, \pi^* \hat\alpha) = \bigl\langle \Dbcoh(\hat X, \hat\alpha), \text{ 10 copies of } \Dbcoh(\P^1, \hat\alpha|_{\P^1}) \bigr\rangle. \]
Now $\hat\alpha|_{\P^1}$ is trivial because $\Br(\P^1) = 0$, and we know that $\Dbcoh(\P^1)$ has a full exceptional collection of length 2, so $\Dbcoh(X, \pi^* \hat\alpha)$ has a full exceptional collection of length $4 + 10 \cdot 2 = 24$.

It remains to show that $\pi^* \hat\alpha = \alpha$.  Because $\Br(X) = \Z/2$, it is enough to show that $\pi^* \hat\alpha$ is not trivial.  If it were trivial, then we would have proved that $\Dbcoh(X)$ has a full exceptional collection, contradicting the fact that $\dimQL(X) = 2$.
\end{proof}

\subsection{Fano 3-folds} \label{fano_3folds}

Iskovskikh classified Fano 3-folds of Picard rank 1 into 17 deformation classes, and Mori and Mukai classified those of higher Picard rank into a further 88 deformation classes.  In every case, $H^3(X,\Z)$ is torsion-free \cite[p.~168]{ip_fano_book}, so by Theorem \ref{3fold_thm} they all satisfy $\dimQL(X) \le 1$.  Beauville's survey paper \cite[\S2.3]{beauville_lueroth_survey} summarizes the situation with rationality in Picard rank 1: in eight of the 17 families, every member is rational; in seven, every member is irrational; and in two, the generic member is irrational.  Belmans' online database \cite{fanography} contains a wealth of information in all Picard ranks.

Kuznetsov has studied derived categories of Fano 3-folds of Picard rank 1 in depth; see \cite[\S2.4]{kuz_rationality_survey} for a survey and references.  In the eight families where every member is rational, he has produced either an exceptional collection of length four, or a semi-orthogonal decomposition consisting of two exceptional bundles and the derived category of a curve.  In six of the other families, he has defined an admissible subcategory $\A_X \subset \Dbcoh(X)$ as the orthogonal to an exceptional collection of length 2; its Hochschild homology looks like that of a curve, although its Hochschild cohomology does not \cite[Thm.~8.9]{kuz_hochschild}, nor does the Serre functor.  In the remaining three families, $\A_X$ is just defined as the orthogonal to the structure sheaf $\O_X$, and its Hochschild homology is like that of a curve except $hh_0 = 3$ rather than 2.

Here we observe that the algebraic K-theory of $\A_X$ also looks like that of a curve.  By Proposition \ref{sod_vs_QL}, $\dimQL(\A_X) \le 1$.  With $\Z$ coefficients, the K-theory of a smooth complex projective curve $C$ of genus $g$ looks like
\[ K_n(C) = \begin{cases}
\Z^2 \oplus J(C) & \text{if $n=0$,} \\
(\Q/\Z)^2 \oplus K_n(C,\Q) & \text{if $n$ is odd and $n \ge 1$, and} \\
(\Q/\Z)^{2g} \oplus K_n(C,\Q) & \text{if $n$ is even and $n \ge 2$,}
\end{cases} \]
where $J(C)$ is the Jacobian of $C$, which is a principally polarized Abelian $g$-fold; for $n=0$ this is well-known, and for $n \ge 1$ it was proved by Pedrini and Weibel \cite[Thm.~3.2]{pw_surfaces}, or see Proposition \ref{pw_prop}.

\begin{prop} \label{3fold_KAX}
With reference to \cite[Table 1]{kuz_rationality_survey}:
\begin{enumerate}
\item If $X$ is one of the Fano 3-folds $X_{14}$, $V_3$, $X_{10}$, $dS_4$, $V_{2,3}$, or $V_6^{1,1,1,2,3}$ and $\A_X \subset \Dbcoh(X)$ is Kuznetsov's subcategory, then
\[ K_n(\A_X) = \begin{cases}
\Z^2 \oplus IJ(X) & \text{if $n=0$,} \\
(\Q/\Z)^2 \oplus K_n(\A_X,\Q) & \text{if $n$ is odd and $n \ge 1$, and} \\
(\Q/\Z)^{2g} \oplus K_n(\A_X,\Q) & \text{if $n$ is even and $n \ge 2$,}
\end{cases} \]
where $g = b_3(X)/2$ and $IJ(X) = H^{1,2}(X)/H^3(X,\Z)$ is the intermediate Jacobian, which is a principally polarized Abelian $g$-fold.

\item If $X$ is $V_{2,2,2}$, $V_4$, or $dS_6$, then
\[ K_n(\A_X) = \begin{cases}
\Z^3 \oplus IJ(X) & \text{if $n=0$,} \\
(\Q/\Z)^3 \oplus K_n(\A_X,\Q) & \text{if $n$ is odd and $n \ge 1$, and} \\
(\Q/\Z)^{2g} \oplus K_n(\A_X,\Q) & \text{if $n$ is even and $n \ge 2$.}
\end{cases} \]
\end{enumerate}
\end{prop}
\begin{proof}
First, the Chow groups of $X$ are given by
\[ CH^i(X) = \begin{cases}
\Z & \text{i=0,} \\
\Z & \text{i=1,} \\
\Z \oplus IJ(X) & \text{i=2,} \\
\Z & \text{i=3.}
\end{cases} \]
For cubic 3-folds (denoted $V_3$ above) this is \cite[\S7.4.1]{huybrechts_cubic_book}, and for the others the proof is similar, using Bloch and Srinivas's \cite[Thm.~1]{bloch-srinivas}. \medskip

Next we claim that
\begin{equation} \label{3fold_K0}
K_0(X) = \Z^4 \oplus IJ(X).
\end{equation}
This follows from Pedrini and Weibel's \cite[Example 6.1.2]{pw_varieties}, but we give the details a little more carefully.

Take the filtration by codimension of support,
\[ 0 = F^4 \subset F^3 \subset F^2 \subset F^1 \subset F^0 = K_0(X). \]
By \cite[Ex.~15.1.5 and 15.3.6]{fulton}, there is a surjection
\[ CH^i(X) \twoheadrightarrow F^i/F^{i+1} \]
which is an isomorphism for $i \le 2$ and whose kernel is 2-torsion for $i=3$.  In this case it is an isomorphism for all $i$, because $CH^3(X)$ is torsion-free.

The inclusion of $F^3 = \Z$ in $K_0(X)$ is split by the Euler characteristic $\chi\colon K_0(X) \to \Z$, so we must have $F^2 = \Z^2 \oplus IJ(X)$.  The rest of the filtration must split as well because the quotients $F^1/F^2$ and $F^0/F^1$ are free Abelian groups.  Thus \eqref{3fold_K0} is established. \medskip

In case (a) we have
\[ K_0(X) = K_0(\A_X) \oplus K_0(\Spec \C)^{\oplus 2} = K_0(\A_X) \oplus \Z^2, \]
so we must have
\[ K_0(\A_X) = \Z^2 \oplus IJ(X) \]
because there are no non-zero maps $IJ(X) \to \Z$.  In case (b) we have $K_0(X) = K_0(\A_X) \oplus K_0(\Spec \C)$, so $K_0(\A_X) = \Z^3 \oplus IJ(X)$ by the same reasoning.

The cohomology of $X$ is torsion-free, as we have seen, so the Atiyah--Hirzebruch spectral sequence degenerates and gives
\begin{align*}
KU^\text{even}(X) &= \Z^4 & KU^\text{odd}(X) &= \Z^{2g}.
\end{align*}
In case (a) we again have $KU^*(X) = KU^*(\A_X) \oplus KU^*(\text{point})^{\oplus 2}$, so
\begin{align*}
KU^\text{even}(\A_X) &= \Z^2 & KU^\text{odd}(\A_X) &= \Z^{2g}
\end{align*}
and in case (b) we get $KU^\text{even}(\A_X) = \Z^3$.  Thus $K_{\ge 1}(\A_X)$ follows from Proposition \ref{pw_prop}.
\end{proof}

\subsection{Enriques categories}

For an Enriques surface $X$, the canonical bundle $\omega_X$ satisfies $\omega_X^{\otimes 2} = \O_X$ but $\omega_X \ne \O_X$; thus on the derived category $\Dbcoh(X)$, the Serre functor $S_X = \underline\quad \otimes \omega_X[2]$ satisfies $S_X^2 = [4]$ but $S_X \ne [2]$.  The Kuznetsov subcategories for some Fano 3-folds enjoy the same property: the main examples \cite[\S3.3]{ppz_enriques} are quartic double solids, which were denoted $dS_4$ earlier; Gushel--Mukai 3-folds, which were denoted $X_{10}$; and Verra 3-folds, which are typically divisors of bidegree $(2,2)$ in $\P^2 \times \P^2$, so they have Picard rank 2, but like the first two they are irrational \cite[\S12.3]{ip_fano_book}.  Perry, Pertusi, and Zhao \cite{ppz_enriques} have studied moduli spaces of Bridgeland-stable objects in these categories, which behave like moduli spaces of Bridgeland-stable objects on Enriques surfaces.

We remark that in each case, the Enriques category $\A_X \subset \Dbcoh(X)$ has $\dimQL(\A_X) = 1$ by Theorem \ref{3fold_thm} and Proposition \ref{sod_vs_QL}, in contrast to the derived category of an actual Enriques surface where $\dimQL(S) = 2$ as we saw in \S\ref{enriques}.  On the other hand, the derived category of twisted sheaves on an Enriques surface is an Enriques category, and we saw in Proposition \ref{twisted_enriques} that $\dimQL(S,\alpha) = 0$.  We are not aware of any work on moduli spaces of twisted sheaves on Enriques surfaces; perhaps they are similar to moduli spaces of untwisted sheaves.

%% file: 4folds.tex

\section{Rationally connected 4-folds} \label{4-folds}

Now we come to our main results.

\subsection{Analysis of the spectral sequence} \label{4fold_thm_sec}

{ \renewcommand{\thethm}{\ref{4fold_thm}}
\addtocounter{thm}{-1}
\begin{thm}
Let $X$ be a smooth complex projective 4-fold that is rationally connected.
\begin{enumerate}
\item $\dimQL(X) \le 3$ if and only if algebraic and homological equivalence coincide on $CH_1(X)$, and in the coniveau filtration
\[ N^2 H^5(X,\Z) \subset N^1 H^5(X,\Z) = H^5(X,\Z), \]
the inclusion is an equality.

\item $\dimQL(X) \le 2$ if and only if in addition, the integral Hodge conjecture holds for $H^4(X,\Z)$ and $H^6(X,\Z)$.

\item If in addition $h^{1,3}(X) = 0$, and $H^5(X,\Z)$ and $H^3(X,\Z)$ are torsion-free, then $\dimQL(X) \le 1$.

\item[(c$'$)] If $\dimQL(X) \le 1$ then $h^{1,3}(X) = 0$, and $H^5(X,\Z)$ is torsion-free.
\end{enumerate}
\end{thm}
}

\begin{proof} First we argue that the $E_2$ page of the Bloch--Ogus--Leray spectral sequence \eqref{bol_ss} with $\Z$ coefficients looks like this:
\[ \xymatrix@=1ex{
4 & 0 & H^1(\H^4_X(\Z)) \ar[rrd] & \dfrac{H^6(X,\Z)}{\text{\footnotesize alg.\,classes}} & 0 & \Z \\
3 & 0 & \dfrac{H^4(X,\Z)}{\text{\footnotesize alg.\,classes}} & N^2 H^5(X,\Z) & \dfrac{CH^3(X)}{\text{\footnotesize alg.\,equiv.}} \\
2 & 0 & H^3(X,\Z) & \dfrac{CH^2(X)}{\text{\footnotesize alg.\,equiv.}} \\
1 & 0 & \NS(X) \\
q=0 & \Z \\
\ar@{-}[]+<1.6em,.8em>;[rrrrr]+<1.6em,.8em> 
\ar@{-}[]+<1.6em,.8em>;[uuuuu]+<1.6em,.8em> 
& p=0 & 1 & 2 & 3 & 4
} \]

As in the proof of Theorem \ref{3fold_thm}, we have $H^0(\H^q_X(\Z)) = 0$ for $q > 0$ by \cite[Prop.~3.3(i)]{ctv}.

The diagonal $p+q=4$ again gives a 2-step filtration of $H^4(X,\Z)$ where the sub-object is $CH^2(X)$ modulo algebraic equivalence, so the quotient $H^1(\H^3_X(\Z))$ is the cokernel of the cycle class map $CH^2(X) \to H^4(X,\Z)$, but now this cycle class map need not be surjective.

The diagonal $p+q=6$ is similar; the differential mapping to $H^3(\H^3_X(\Z))$ might not be zero, but the quotient of the 2-step filtration on $H^6(X,\Z)$ is still the cokernel of the cycle class map $CH^3(X) \to H^6(X,\Z)$.

On the diagonal $p+q=7$, we again have $H^7(X,\Z) = 0$ because $X$ is simply connected.

Thus we have justified the $E_2$ page shown above. \medskip

(a) From \eqref{K/tau_d} we have
\[ (K/\tau)_4(X,\Z/m) = H^0(\H^4_X(\Z/m)). \]
The group $H^1(\H^4_X(\Z))$ is torsion by \cite[Prop.~3.3(ii)]{ctv}, so from the long exact sequence \eqref{ctv_les} we see that $H^0(\H^4_X(\Z/m)) = 0$ for all $m$ if and only if $H^1(\H^4_X(\Z)) = 0$.

Looking at the diagonals $p+q=5$ and $p+q=6$ in our $E_2$ page above, we get an exact sequence
\begin{equation} \label{H1H4seq}
0 \to N^2 H^5(X,\Z) \to H^5(X,\Z) \to H^1(\H^4_X(\Z)) \to \frac{CH^3(X)}{\text{\footnotesize alg.\,equiv.}} \to H^6(X,\Z),
\end{equation}
so $H^1(\H^4_X(\Z)) = 0$ if and only if the first map is an isomorphism and the last is injective; but the kernel of the last map is exactly the quotient of homologically trivial cycles by algebraically trivial ones.

Thus we have proved that $(K/\tau)_4(X,\Z/m) = 0$ for all $m$ if and only if algebraic and homological equivalence coincide on $CH^3(X)$, and $N^2 H^5(X,\Z) = H^5(X,\Z)$.\medskip

(b) We have $H^3(\H^4_X(\Z)) = 0$ and now $H^1(\H^4_X(\Z)) = 0$.  We know that $H^2(\H^4_X(\Z))$ is finitely generated because it is a quotient of $H^6(X,\Z)$, and torsion because the Hodge conjecture holds for $H^6$ of a 4-fold.  Thus we see from the long exact sequence \eqref{ctv_les} that the following are equivalent:
\begin{itemize}
\item $H^1(\H^4_X(\Z/m)) = 0$ for all $m$,
\item $H^2(\H^4_X(\Z)) = 0$,
\item $H^2(\H^4_X(\Z/m)) = 0$ for all $m$.
\end{itemize}
The exact sequence \eqref{K/tau_d-1} reads
\begin{multline*}
0 \to H^1(\H^4_X(\Z/m)) \to (K/\tau)_3(X, \Z/m) \\
\to H^0(\H^3_X(\Z/m)) \to H^2(\H^4_X(\Z/m)),
\end{multline*}
so we see that $(K/\tau)_3(X, \Z/m) = 0$ for all $m$ if and only if $H^2(\H^4_X(\Z)) = 0$ and $H^0(\H^3_X(\Z/m)) = 0$ for all $m$.  The first is equivalent to saying that the cycle class map $CH^3(X) \to H^6(X,\Z)$ is surjective.  The second is equivalent to saying that $H^1(\H^3_X(\Z))$ is torsion-free, or that the cokernel of the cycle class map $CH^2(X) \to H^4(X,\Z)$ is torsion-free.  The (rational) Hodge conjecture holds for rationally connected (or just uniruled) 4-folds \cite{conte-murre}, so this last condition is exactly the integral Hodge conjecture. \medskip

(c) and (c$'$) Now the long exact sequence \eqref{ctv_les} gives
\[ H^0(\H^2_X(\Z/m)) = m\text{-tors}(H^3(X,\Z)) \]
and
\[ 0 \to \dfrac{H^4(X,\Z)}{\text{\footnotesize alg.\,classes}} \otimes \Z/m \to H^1(\H^3_X(\Z/m)) \to m\text{-tors}(H^5(X,\Z)) \to 0. \]
The first term vanishes for all $m$ if and only if $H^4(X,\Z)$ is generated by algebraic classes, and we have seen that the integral Hodge conjecture holds, so this is true if and only if $h^{1,3} = 0$.

From the $E_2$ page of \eqref{our_E_2} and the vanishing of various $H^i(\H^j_X(\Z/m))$ established earlier, we get an exact sequence
\begin{multline*}
0 \to H^1(\H^3_X(\Z/m)) \to (K/\tau)_2(X, \Z/m) \\
\to H^0(\H^2_X(\Z/m)) \to H^2(\H^3_X(\Z/m)).
\end{multline*}
Now (c) amounts to saying that if the first and third terms vanish then the second does, and (c$'$) amounts to saying that if the second term vanishes then the first does.
\end{proof}

\subsection{Cubic and Gushel--Mukai 4-folds} \label{cubic_gm_sec}

Many smooth cubic hypersurfaces $X \subset \P^5$ are known to be rational, and a very general one is expected to be irrational, although none has yet been proven to be irrational.\footnote{While this paper was under review, the long-awaited \cite{kkpy} was posted to the arXiv, using quantum cohomology to say that a very general cubic 4-fold is irrational.  The argument has been refined and streamlined in \cite{guere}.}  Kuznetsov \cite{kuz_cubics} introduced a semi-orthogonal decomposition
\[ \Dbcoh(X) = \langle \A_X, \O_X, \O_X(1), \O_X(2) \rangle, \]
where $\A_X$ behaves like a non-commutative K3 surface; it has been heavily studied in connection with hyperk\"ahler geometry as well as classical rationality questions.

A Gushel--Mukai 4-fold is either a smooth intersection of the Grassmannian $\Gr(2,5) \subset \P^9$ with a hyperplane and a quadric, or the double cover of the intersection of two hyperplanes in $\Gr(2,5)$ branched over its intersection with a quadric.  The interest in GM 4-folds stems from their remarkable similarity to cubic 4-folds, with connections to K3 surfaces and hyperk\"ahler varieties through their cohomology, derived categories of coherent sheaves, and various geometric constructions \cite{debarre_gm_survey}.  The situation with rationality is also the same: some GM 4-folds are known to be rational, and the very general one is expected to be irrational, but none has yet been proven to be irrational.\footnote{The techniques of \cite{kkpy} and \cite{guere} have been applied to Gushel--Mukai 4-folds in \cite{bmp}.}  There are also GM 4-folds that are expected to be irrational and are birational to cubic 4-folds \cite[Thm.~5.8]{kuz_perry}.  The Kuznetsov subcategory for a GM 4-fold is defined by a semi-orthogonal decomposition
\[ \Dbcoh(X) = \langle \A_X, U_X, \O_X, U_X(1), \O_X(1) \rangle, \]
where $U_X$ is the restriction (or pullback) of the universal rank 2 bundle from $\Gr(2,5)$.

\begin{prop} \label{cubics_are_QL2}
Every smooth complex cubic 4-fold or Gushel--Mukai 4-fold $X$ has $\dimQL(X) = 2$.
\end{prop}
\begin{proof}
We have $h^{1,3}(X) = 1$, so $\dimQL(X) \ge 2$ by Theorem \ref{4fold_thm}(c$'$).

To see that $\dimQL(X) \le 3$, we can apply Proposition \ref{conic_fibr} and Corollary \ref{K/tau_birational}: a cubic 4-fold is birational to a conic fibration by projecting from a line in $X$ onto a complementary $\P^3$, and a GM 4-fold is birational to a conic fibration by \cite[\S3]{dim} or \cite[\S2.1]{ppz}.\footnote{We thank Lie Fu for this argument.}

We give another proof that a cubic 4-fold has $\dimQL(X) \le 3$, using Theorem \ref{4fold_thm}(a).  Shen proved that $CH_1(X)$ is generated by lines \cite[Thm.~1.1]{shen_lines}, and Tian and Zong gave another proof \cite[Thm.~1.7]{tian_zong}, so if $F$ denotes the variety of lines on $X$, then the cylinder map $CH_0(F) \to CH_1(X)$ is surjective; homological and algebraic equivalence coincide on $CH_0(F)$, so they coincide on $CH_1(X)$ as well.  And we have $H^5(X,\Z) = 0$ by the Lefschetz hyperplane theorem.

Now to see that $\dimQL(X) \le 2$ in both cases, we use Theorem \ref{4fold_thm}(b).  The integral Hodge conjecture for $H^4$ of cubic 4-folds was proved by Voisin \cite[Thm.~18]{voisin_Z_hodge}, and later re-proved by Mongardi and Ottem \cite[Cor.~0.3]{mongardi-ottem}, and by Perry who also proved it for GM 4-folds \cite[Cor.~1.2]{perry}.  And in both cases $H^6(X,\Z) = \Z$ is generated by the class of a line.  
\end{proof}

Thus it seems that higher K-theory gives no obstruction to the rationality of a general cubic 4-fold, which was the original motivation for this paper.  But we can recast this birational disappointment into a positive statement, and describe the higher K-theory of Kuznetsov's K3 category $\A_X$ as completely as we can that of an honest K3 surface.

If $S$ is a (complex) K3 surface, then
\[ K_n(S) = \begin{cases}
\Z^{r+2} \oplus CH_0(S)_\text{hom} & \text{if $n = 0$,} \\
(\Q/\Z)^{24} \oplus K_n(S,\Q) & \text{ if $n$ is odd and $n \ge 1$, and} \\
K_n(S,\Q) & \text{ if $n$ is even and $n \ge 2$,}
\end{cases} \]
where $r = \rank(\NS(S))$ and $CH_0(S)_\text{hom} = CH_0(S)_\text{alg}$ is a $\Q$-vector space by Roitman's theorem.  For $n=0$ this is \cite[Cor.~12.1.5]{huybrechts_k3}, and for $n \ge 1$ it follows from Pedrini and Weibel's \cite[Example~6.7]{pw_surfaces}, or see Proposition \ref{pw_prop}.  The K-theory of $\A_X$ is entirely similar:

{ \renewcommand{\thethm}{\ref{4fold_KAX}(b)}
\addtocounter{thm}{-1}
\begin{thm}
Let $X$ be a smooth complex cubic 4-fold, and let $\A_X \subset \Dbcoh(X)$ be Kuznetsov's K3 category.  Then
\[ K_n(\A_X) = \begin{cases}
\Z^{r+1} \oplus CH_1(X)_\mathrm{hom} & \text{if $n = 0$,} \\
(\Q/\Z)^{24} \oplus K_n(\A_X,\Q) & \text{ if $n$ is odd and $n \ge 1$, and} \\
K_n(\A_X,\Q) & \text{ if $n$ is even and $n \ge 2$,}
\end{cases} \]
where $r = \rank(H^4(X,\Z) \cap H^{2,2}(X))$ and $CH_1(X)_\mathrm{hom} = CH_1(X)_\mathrm{alg}$ is a $\Q$-vector space.

The same holds for a GM 4-fold with $\Z^r$ in place of $\Z^{r+1}$.
\end{thm}
}
\begin{proof}
The Chow groups of a cubic 4-fold $X$ are given in \cite[\S7.4.1]{huybrechts_cubic_book}:
\[ CH^i(X) = \begin{cases}
\Z & i = 0, \\
\Z & i = 1, \\
\Z^r & i = 2, \\
\Z \oplus CH_1(X)_\text{hom} & i = 3, \\
\Z & i = 4.
\end{cases} \]
We saw in the proof of Proposition \ref{cubics_are_QL2} that $CH_1(X)_\text{hom} = CH_1(X)_\text{alg}$, so it is a divisible group.  One can deduce that it is torsion-free from a careful geometric analysis of the kernel of the cylinder map (see \cite[Rmk.~7.4.1(ii)]{huybrechts_cubic_book}), or more simply from Proposition \ref{CH1} below.

The Chow groups of a GM 4-fold admit an identical description: we can emulate the computation in \cite[\S7.4.1]{huybrechts_cubic_book}, saying that $\Pic(X) = \Z$ by the Lefschetz hyperplane theorem, that the cycle map $CH^2(X) \to H^4(X,\Z) \cap H^{2,2}(X)$ is injective by Bloch--Srinivas \cite{bloch-srinivas} and surjective by Perry's verification of the integral Hodge conjecture \cite[Cor.~1.2]{perry}, that $CH_1(X)_\text{alg} = CH_1(X)_\text{hom}$ by Proposition \ref{cubics_are_QL2}, and that this is torsion-free by Proposition \ref{CH1} below.  Recall that $CH^k(X)_\text{alg}$ is always divisible \cite[Lec.~1, Lem.~1.3]{bloch_book}. \medskip

Next we compute $K_0(X)$.  As in the proof of Proposition \ref{3fold_KAX}, we have the filtration by codimension of support,
\[ 0 = F^5 \subset F^4 \subset F^3 \subset F^2 \subset F^1 \subset F^0 = K_0(X), \]
and a surjection
\[ CH^i(X) \twoheadrightarrow F^i/F^{i+1} \]
which is an isomorphism for $i \le 2$ and whose kernel is 2-torsion for $i=3$ and 6-torsion for $i=4$, by \cite[Ex.~15.3.6]{fulton}.  Thus it is an isomorphism for all $i$ in this case, because we have seen that $CH^*(X)$ is torsion-free.  The inclusion of $F_4 = \Z$ in $K_0(X)$ is split by the Euler characteristic $\chi\colon K_0(X) \to \Z$, and the rest of the filtration must be split because the quotients are free Abelian groups at every step.  Thus we get
\[ K_0(X) = \Z^{r+4} \oplus CH_1(X)_\text{hom}. \]

For a cubic 4-fold, we have
\[ K_0(X) = K_0(\A_X) \oplus K_0(\Spec \C)^{\oplus 3}, \]
so as in the proof of Proposition \ref{3fold_KAX} we must have
\[ K_0(\A_X) = \Z^{r+1} \oplus CH_1(X)_\text{hom} \]
because there are no non-zero maps from the divisible group $CH_1(X)_\text{hom}$ to $\Z$.  For a GM 4-fold it is the same but with $K_0(\Spec \C)^{\oplus 4}$.

For the higher K-theory, we have $KU^\text{even}(X) = \Z^{27}$ for a cubic 4-fold or $\Z^{28}$ for a GM 4-fold, and $KU^\text{odd}(X) = 0$ in both cases, so $KU^\text{even}(\A_X) = \Z^{24}$ and $KU^\text{odd}(\A_X) = 0$.  Thus $K_{\ge 1}(\A_X)$ follows from Proposition \ref{pw_prop}.
\end{proof}

\begin{prop} \label{CH1}
Suppose that $X$ is a smooth complex projective variety with $CH_0(X) = \Z$.
If $H^5(X,\Z) = 0$, then $CH^3(X)_\mathrm{alg}$ is torsion-free.  If $H_3(X,\Z) = 0$, then $CH_1(X)_\text{alg}$ is torsion-free.
\end{prop}
\begin{proof}
Let $d = \dim(X)$.  By Poincar\'e duality $H^{2d-3}(X,\Z) = H_3(X,\Z)$, and by definition $CH^{d-1}(X) = CH_1(X)$.

Ma proved in \cite[Thm.~5.1]{ma} that the torsion in $CH^p(X)_\text{alg}$ is a quotient of $H^{p-1}(\H^p_X(\Z)) \otimes \Q/\Z$.  The Bloch--Ogus--Leray spectral sequence yields exact sequences
\begin{gather*}
H^0(\H^4_X(\Z)) \to H^2(\H^3_X(\Z)) \to H^5(X,\Z) \\
H^{d-4}(\H^d_X(\Z)) \to H^{d-2}(\H^{d-1}_X(\Z)) \to H^{2d-3}(X,\Z).
\end{gather*}
If $CH_0(X) = \Z$ then $H^0(\H^4_X(\Z)) = 0$ by \cite[Prop.~3.3(i)]{ctv}, and if $H^5(X,\Z) = 0$ then $H^2(\H^3_X(\Z)) = 0$, which proves the first claim.  If $H^{2d-3}(X,\Z) = 0$ then $H^{d-2}(\H^{d-1}_X(\Z))$ is a quotient of $H^{d-4}(\H^d_X(\Z))$, which is torsion by \cite[Prop.~3.3(ii)]{ctv}, so after tensoring with the divisible group $\Q/\Z$ it becomes zero, which proves the second claim.

We give another proof, due to Claire Voisin. In \cite[Thm.~2.17]{voisin_birational} she proved that if a variety $X$ admits a Chow decomposition of the diagonal
\[ \Delta_X = X \times x + Z \text{ in } CH^d(X \times X), \]
where $x$ is any point and $Z$ is supported in $D \times X$ for some divisor $D \subset X$, then the kernels of the Abel--Jacobi maps
\begin{gather*}
CH^3(X)_\text{hom} \to J^5(X) := H^5(X,\C)/(F^3 + H^5(X,\Z)_\text{tf}) \\
CH^{d-1}(X)_\text{hom} \to J^{2d-3}(X) := H^{2d-3}(X,\C)/(F^{d-1} + H^{2d-3}(X,\Z)_\text{tf})
\end{gather*}
are torsion-free.  The same proof shows that if $X$ admits a \emph{rational} decomposition of the diagonal, meaning that
\[ N \Delta_X = N(X \times x) + Z \text{ in } CH^d(X \times X) \]
for some integer $N$, then the torsion in these kernels is annihilated by $N$.  We have $CH_0(X) = \Z$, so $X$ admits a rational decomposition of the diagonal by \cite[Prop.~1]{bloch-srinivas}.  If $H^5(X,\Z) = 0$ then $J^5(X) = 0$, so the torsion in $CH^3(X)_\text{hom}$ has bounded order; thus the same is true of the subgroup $CH^3(X)_\text{alg}$, but this is a divisible group, so its torsion must vanish.  Similarly, if $H^{2d-3}(X,\Z) = 0$ then $CH^{d-1}(X)_\text{alg}$ is torsion-free.
\end{proof}

\subsection{Quartic 4-folds}

Totaro showed that a very general quartic 4-fold is irrational \cite[Thm.~2.1]{totaro_irr}, and even gave an explicit example of an irrational one defined over $\Q$ \cite[Example 3.1]{totaro_irr}.  We do not know whether any smooth quartic 4-fold is rational.

We use Theorem \ref{4fold_thm} to see that every quartic 4-fold $X$ satisfies $2 \le \dimQL(X) \le 3$, and a very general one satisfies $\dimQL(X) = 2$:
\begin{itemize}
\item $CH_1(X)$ is generated by lines by \cite[Thm.~1.7]{tian_zong}, so homological and algebraic equivalence coincide by the same cylinder map argument as in the proof of Proposition \ref{cubics_are_QL2}.
\item $H^5(X,\Z) = 0$ by the Lefschetz hyperplane theorem.
\item $H^6(X,\Z) = \Z$ is generated by the class of a line.
\item If $X$ is very general then $H^4(X,\Z) \cap H^{2,2}(X) = \Z \cdot h^2$, where $h \in H^2(X,\Z)$ is the hyperplane class.  But if $X$ is ``Noether--Lefschetz special'' then we do not know whether the integral Hodge conjecture holds on $H^4(X,\Z)$.
\item $h^{1,3}(X) = 21$.
\end{itemize}

\subsection{Hassett, Pirutka, and Tschinkel's 4-folds}

In \cite{hpt1}, Hassett, Pirutka, and Tschinkel proved that a very general hypersurface of bi\-degree $(2,2)$ in $\P^2 \times \P^3$ is irrational, but many are rational, and the rational ones form a dense subset of the moduli space (in the analytic topology): in short, they exhibit the behavior expected of cubic 4-folds.  In \cite{hpt2}, they proved the same for complete intersections of 3 quadrics in $\P^7$.  We verify that any smooth 4-fold $X$ in either family has $\dimQL(X) = 2$.

On the one hand, we can check the hypotheses of Theorem \ref{4fold_thm}.  Algebraic and homological equivalence coincide on $CH_1$ as in the proof of Proposition \ref{cubics_are_QL2}: a $(2,2)$ divisor in $\P^2 \times \P^3$ admits a conic fibration by projecting onto $\P^3$, and an intersection of three quadrics in $\P^7$ has $CH_1$ generated by lines by Tian and Zong's \cite[Thm.~1.7]{tian_zong}.  In both cases we have $H^5(X,\Z) = 0$ by the Lefschetz hyperplane theorem, and $H^6(X,\Z)$ is generated by algebraic classes, either by hand or by \cite[Thm.~1.7(ii)]{voisin-horing}.  The integral Hodge conjecture for $H^4(X,\Z)$ is discussed in \cite[Rmk.~7]{hpt1} and \cite[\S2.3]{hpt2}.

Alternatively, we can see that $\dimQL(X) = 2$ using Proposition \ref{sod_vs_QL}, Corollary \ref{twisted}, and Kuznetsov's work on derived categories of quadric fibrations and intersections of quadrics \cite{kuznetsov_quadrics}, as follows.  For a complete intersection of three quadrics in $\P^7$, by \cite[Thm.~5.5 and Prop.~3.13]{kuznetsov_quadrics} there is a semi-orthogonal decomposition
\[ \Dbcoh(X) = \langle \Dbcoh(S,\alpha), \O_X, \O_X(1) \rangle, \]
where $S$ is the double cover of $\P^2$ branched over a smooth octic curve and $\alpha$ is a Brauer class of index 2; see also \cite{nick_thesis}.  For a $(2,2)$ divisor in $\P^2 \times \P^3$, the projection onto $\P^2$ makes $X$ into a quadric surface fibration degenerating over an octic curve, and by \cite[Thm.~4.2 and Prop.~3.13]{kuznetsov_quadrics} there is a semi-orthogonal decomposition 
\begin{multline*}
\Dbcoh(X) = \langle \Dbcoh(S,\alpha),\ \O_X(0,1), \O_X(1,1), \O_X(2,1), \\
\O_X(0,2), \O_X(1,2), \O_X(2,2) \rangle,
\end{multline*}
where $S$ is the double cover of $\P^2$ branched over the octic curve and $\alpha$ is again a Brauer class of index 2.
In either case, Proposition \ref{sod_vs_QL} gives $\dimQL(X) = \dimQL(S,\alpha)$, and Corollary \ref{twisted} gives $\dimQL(S,\alpha) \le 2$.

\subsection{Schreieder's 4-fold}

Having seen many rationally connected 4-folds with $\dimQL(X) = 2$, one begins to wonder whether $\dimQL(X) > 2$ is possible.  In \cite[Cor.~1.6]{schreieder_small_slopes}, Schreieder constructed a (smooth complex projective) unirational 4-fold $X$ that does not satisfy the integral Hodge conjecture on $H^4(X,\Z)$; thus $\dimQL(X) \ge 3$ by Theorem \ref{4fold_thm}(b).  It admits a map $X \to \P^3$ whose generic fiber is a conic, so Proposition \ref{conic_fibr} gives $\dimQL(X) \le 3$ as well.

Schreieder has also constructed higher-dimensional examples with non-zero unramified cohomology in various degrees \cite[Thm.~1.5]{schreieder_small_slopes}, but it is unclear how this non-vanishing filters through our spectral sequence to affect the Lichtenbaum--Quillen dimension of the varieties in question.

\subsection{Homological projective duality examples}

\newcommand{\fanoFourFoldHodgeDiamond}[3]{
  \[ 
  \scriptsize
  \setlength \arraycolsep {1.5pt}
  \renewcommand \arraystretch {.8}
  \begin{array}{ccccccccc}
  & & & & 1 \\
  & & & 0 & & 0 \\
  & & 0 & & 1 & & 0 \\
  & 0 & & 0 & & 0 & & 0 \\
  0 & & #1 & & #2 & & #3 & & 0. \\
  & 0 & & 0 & & 0 & & 0 \\
  & & 0 & & 1 & & 0 \\
  & & & 0 & & 0 \\
  & & & & 1
  \end{array}
  \]
}

\newcommand{\surfaceHodgeDiamond}[3]{
  \[
  \scriptsize
  \setlength \arraycolsep {1.5pt}
  \renewcommand \arraystretch {.8}
  \begin{array}{ccccc}
  & & 1 \\
  & 0 & & 0 \\
  #1 & & #2 & & #3. \\
  & 0 & & 0 \\
  & & 1
  \end{array}
  \]
}

In this section we prove:
{ \renewcommand{\thethm}{\ref{hpd_thm}}
\addtocounter{thm}{-1}
\begin{thm}
The following Fano 4-folds $X$ satisfy the integral Hodge conjecture on $H^4(X,\Z)$, and $\Griff_1(X) = 0$.
\begin{enumerate}
\item Intersections of $\Gr(2,7)$ with 6 general hyperplanes in the Pl\"ucker embedding.
\item ``Pfaffian'' 4-folds obtained as linear sections of the space of $7 \times 7$ skew-symmetric matrices of rank 4.
\item The linear sections of the double quintic symmetroid studied by Ottem and Rennemo in \cite{ottem-rennemo}.
\end{enumerate}
\end{thm}
}
\noindent In each case we will describe the Fano 4-fold $X$ in detail, and an associated surface $Y$ of general type with $h^{0,2}(Y) \ne 0$.  Kuznetsov's homological projective duality will give a semi-orthogonal decomposition of $\Dbcoh(X)$ into a copy of $\Dbcoh(Y)$ and three or four exceptional objects.  Then Proposition \ref{sod_vs_QL} and Theorem \ref{surface_thm}(a) imply that $\dimQL(X) = 2$, so Theorem \ref{4fold_thm}(a) implies that $\Griff_1(X) = 0$, and Theorem \ref{4fold_thm}(b) implies that the integral Hodge conjecture holds for $H^4(X,\Z)$.  (The statements from Theorem \ref{4fold_thm} about $H^5$ and $H^6$ are not interesting in these examples.)

Homological projective duality gives many other examples of 4-folds $X$ with semi-orthogonal decompositions that make $\dimQL(X) \le 2$, but the others that we found were already known to be rational.

\subsubsection{Linear sections of Gr(2,7)}
The Pl\"ucker embedding maps the 10-dimensional $\Gr(2,7)$ into the $\P^{20}$ of $7 \times 7$ skew-symmetric matrices, identifying it with the set of matrices of rank 2.  Intersect $\Gr(2,7)$ with 6 general hyperplanes to get a smooth Fano 4-fold $X$.  We can compute its Hodge diamond using the Lefschetz hyperplane theorem and some Hirzebruch--Riemann--Roch calculations using the Schubert2 package of Macaulay2 \cite{M2}:
\fanoFourFoldHodgeDiamond{6}{57}{6}
These Fano 4-folds are type (b7) in K\"uchle's classification \cite[\S3]{kuechle}.

The classical projective dual of $\Gr(2,7)$ is the 17-dimensional space of matrices of rank $\le 4$, which we denote $\Pf(4,7)$; its singular locus is exactly $\Gr(2,7)$.  Really we should be careful to distinguish between the $\P^{20}$ of skew-symmetric maps $V \to V^*$, where $V \cong \C^7$, and the dual $\P^{20}$ of skew-symmetric maps $V^* \to V$.  Our 6 hyperplanes in the first $\P^{20}$ correspond to 6 points in the dual $\P^{20}$, which span a $\P^5$ that intersects $\Pf(4,7)$ transversely in a smooth surface $Y$ of general type.  Its Hodge diamond is
\surfaceHodgeDiamond{6}{56}{6}

Kuznetsov's \cite[Thm.~4.1]{kuz_gr} gives a semi-orthogonal decomposition
\[ \Dbcoh(X) = \langle \Dbcoh(Y), \Sym^2 U_X, U_X, \O_X \rangle \]
where $U_X$ is the restriction of the universal rank 2 bundle from $\Gr(2,7)$.

\subsubsection{Pfaffian 4-folds}

Retain the set-up of the previous example, but now reverse the roles: let $Y$ be the intersection of $\Gr(2,7)$ with 8 general hyperplanes, which is a smooth surface of general type.  Its Hodge diamond is
\surfaceHodgeDiamond{13}{98}{13}
The 8 hyperplanes in $\P^{20}$ correspond to 8 points in the dual $\P^{20}$, which span a $\P^7$ that intersects $\Pf(4,7)$ transversely in a smooth Fano 4-fold $X$.  Its Hodge diamond is
\fanoFourFoldHodgeDiamond{13}{99}{13}
Now \cite[Thm.~4.1]{kuz_gr} gives a semi-orthogonal decomposition
\[ \Dbcoh(X) = \langle \Dbcoh(Y), F_0^*|_X, F_1^*|_X, \O_X \rangle \]
where $F_0$ and $F_1$ are certain homogeneous bundles on $\Pf(4,7) \setminus \Gr(2,7)$. \medskip

For another description of $X$, we can take the map $\Pf(4,7) \setminus \Gr(2,7) \dashrightarrow \Gr(3,7)$ that sends a rank-4 matrix to its kernel; this embeds $X$ into $\Gr(3,7)$ as the first degeneracy locus of a general map $\O^{13} \to \Lambda^2 Q$, where $Q$ is the universal rank-4 quotient bundle.

\subsubsection{Ottem and Rennemo's 4-folds}

The $\P^{14}$ of $5 \times 5$ symmetric matrices is also stratified by rank, with the dimension and degree of each stratum as shown:\footnote{Harris's book \cite{joe_harris} has a friendly discussion of this stratification in Example 22.31, although he incorrectly gives the degree of the rank $\le 3$ locus as 10, and in Prop.~22.32 the $n-k-\alpha$ should be $n-k-\alpha+1$.}
\[ \begin{array}{ccccccccc}
\{ \rank = 1 \} & \subset & \{ \rank \le 2 \} & \subset & \{ \rank \le 3 \} & \subset & \{ \rank \le 4 \} & \subset & \P^{14} \\
\dim = 4 & & \dim = 8 & & \dim = 11 & & \dim = 13 \\
\deg = 16 & & \deg = 35 & & \deg = 20 & & \deg = 5
\end{array} \]
The rank 1 locus is identified with $\P^4$ in its second Veronese embedding.  The rank $\le 2$ locus is identified with $\Sym^2 \P^4$.  The rank $\le 4$ locus is called the ``quintic symmetroid.''  Each stratum is smooth away from the one before.

The classical projective dual of $\Sym^2 \P^4$ is the rank $\le 3$ locus, but its homological projective dual is the double cover of the quintic symmetroid branched over the rank $\le 3$ locus; to be precise we need a non-commutative resolution of singularities on both sides, but we will take linear sections that avoid the singularities.  Hosono and Takagi studied this example in a series of papers including \cite{ht_towards}, and Rennemo developed it further in \cite{joergen_thesis}.

First we can intersect $\Sym^2 \P^4$ with 6 general hyperplanes to get a smooth surface $Y$ of general type with $\pi_1 = \Z/2$.  Its Hodge diamond is
\surfaceHodgeDiamond{9}{65}{9}
Dually, we can intersect the quintic symmetroid with a $\P^5$ and take its double cover branched over the rank $\le 3$ locus to get a smooth Fano 4-fold $X$ whose Hodge diamond is
\fanoFourFoldHodgeDiamond{9}{67}{9}
Ottem and Rennemo have studied this 4-fold in \cite{ottem-rennemo}, proving that \linebreak 
$H^3(X,\Z) = \Z/2$.  By \cite[Prop.~4.3]{ottem-rennemo}, there is a semi-orthogonal decomposition
\[ \Dbcoh(X) = \langle \Dbcoh(Y), E_1, E_2, E_3, E_4 \rangle, \]
where $E_1, \dotsc, E_4$ are exceptional objects.

%% file: 5folds.tex

\section{Cubic 5-folds} \label{cubic_5fold_sec}

Nothing seems to be known about rationality of smooth cubic 5-folds; it is reasonable to guess that most or all are irrational, but Lichtenbaum--Quillen dimension sheds no light on the question:

{ \renewcommand{\thethm}{\ref{cubic_5fold}}
\addtocounter{thm}{-1}
\begin{thm}
If $X$ is a smooth complex cubic 5-fold, then $\dimQL(X) = 1$.
\end{thm}
}
\begin{proof}
Fu and Tian \cite[Thm.~3]{fu_tian} showed that the $E_2$ page of the Bloch--Ogus--Leray spectral sequence \eqref{bol_ss} with $\Z$ coefficients looks like this:
\[ \xymatrix@=1ex{
5 & 0 & 0 & 0 & 0 & 0 & \Z \\
4 & 0 & 0 & 0 & 0 & \Z \\
3 & 0 & 0 & \Z^{42} & \Z \\
2 & 0 & 0 & \Z \\
1 & 0 & \Z \\
q=0 & \Z \\
\ar@{-}[]+<1.6em,.8em>;[rrrrrr]+<1.6em,.8em> 
\ar@{-}[]+<1.6em,.8em>;[uuuuuu]+<1.6em,.8em> 
& p=0 & 1 & 2 & 3 & \quad4\quad & \quad5\quad
} \]
We immediately compute $H^p(\H^q_X(\Z/m))$ using the long exact sequence \eqref{ctv_les}, and our spectral sequence gives $(K/\tau)_n(X,\Z/m) = 0$ for all $n > 1$.
\end{proof}

The Kuznetsov component of the derived category of a cubic 5-fold is defined by a semi-orthogonal decomposition
\[ \Dbcoh(X) = \langle \A_X, \O_X, \O_X(1), \O_X(2), \O_X(3) \rangle. \]
Its Hochschild homology looks like that of a genus-21 curve, although its Serre functor does not ($S_{\A_X}^3 = [7]$) and its Hochschild cohomology probably does not either.  We conclude by showing that its higher K-theory does look like that of a curve, but the analysis of $K_0$ is more delicate than it was for 3-folds and 4-folds.
\begin{prop}
If $\A_X$ is the Kuznetsov category of a smooth complex cubic 5-fold, then
\[ K_n(\A_X) = \begin{cases}
\Z^2 \oplus IJ(X) & \text{if $n=0$,} \\
(\Q/\Z)^2 \oplus K_n(\A_X,\Q) & \text{if $n$ is odd and $n \ge 1$, and} \\
(\Q/\Z)^{42} \oplus K_n(\A_X,\Q) & \text{if $n$ is even and $n \ge 2$,}
\end{cases} \]
where $IJ(X) = H^{2,3}(X)/H^5(X,\Z)$ is the intermediate Jacobian, which is a principally polarized Abelian 21-fold.
\end{prop}
\begin{proof}
Fu and Tian \cite{fu_tian} have also computed the Chow groups of $X$:
\[ CH^i(X) = \begin{cases}
\Z & \text{i=0,} \\
\Z & \text{i=1,} \\
\Z & \text{i=2,} \\
\Z \oplus IJ(X) & \text{i=3,} \\
\Z & \text{i=4,} \\
\Z & \text{i=5,}
\end{cases} \]
As in the proofs of Proposition \ref{3fold_KAX} and Theorem \ref{4fold_KAX}(b), we have the filtration by codimension of support
\[ 0 = F^6 \subset F^5 \subset F^4 \subset F^3 \subset F^2 \subset F^1 \subset F^0 = K_0(X), \]
and a surjection
\begin{equation} \label{fulton_surjection}
CH^i(X) \twoheadrightarrow F^i/F^{i+1}
\end{equation}
whose kernel is $(i-1)!$-torsion by \cite[Ex.~15.1.5 and 15.3.6]{fulton}.  Thus it is an isomorphism except perhaps for $i=3$, because $CH^3(X)$ has 2-torsion.

To show that it is also an isomorphism for $i=3$, we use the motivic spectral sequence
\[ E_2^{p,q} = H^{p-q}_\text{mot}(X,\Z(-q)) \Longrightarrow K_{-p-q}(X). \]
On the diagonal $p+q=0$ we find $H^{2p}_\text{mot}(X,\Z(p)) = CH^p(X)$ with no differentials going out; the resulting filtration on $K_0(X)$ is the one above, so the surjections \eqref{fulton_surjection} are isomorphisms if and only if the differentials coming in are zero.

On the $E_2$ page, the differential coming into $H^6_\text{mot}(X,\Z(3)) = CH^3(X)$ comes from $H^3_\text{mot}(X,\Z(2))$; let us argue that it is 2-torsion.  It commutes with the Adams operations $\psi^k$ for all $k$, which act on $H^6_\text{mot}(X,\Z(3))[\tfrac1k]$ by multiplication by $k^3$ and on $H^3_\text{mot}(X,\Z(2))[\tfrac1k]$ by multiplication by $k^2$; the most solid reference for our purposes seems to be \cite[Lem.~B.2]{e-morrow}, but see also \cite[Thm.~7]{gillet-soule}.
Thus the differential tensored with $\Z[\tfrac1k]$ is annihilated by $k^3 - k^2 = k^2(k-1)$, hence is annihilated by $k-1$.  Taking $k=2$, we see that the differential tensor $\Z[\tfrac12]$ is zero, so it is 2-power torsion.  Taking $k=3$, we see that it is 2-torsion.

Thus the differential descends to $H^3_\text{mot}(X,\Z(2))/2$, which injects into $H^3_\text{mot}(X,\Z/2(2))$, which injects into $H^3_\text{sing}(X,\Z/2)$ by the Beilinson--Lichtenbaum conjecture, and this vanishes by the Lefschetz hyperplane theorem.  So the relevant differential on the $E_2$ page vanishes as desired.

On the $E_3$ page, the differential coming into $H^6_\text{mot}(X,\Z(3))$ comes from $H^1_\text{mot}(X,\Z(1)) = \Gamma(\O_X^*) = \C^*$, and there are several ways to see that it vanishes.  For one, the differential is 24-torsion (again using Adams operations) and the source is a divisible group.  For another, the pullback map from $K_1(\Spec \C) = \C^*$ to $K_1(X)$ is split by restriction to a point, so this $\C^*$ must survive to the end of the spectral sequence.

After this there are no more differentials coming into $H^6_\text{mot}(X,\Z(3))$.

Next we argue that the codimension filtration on $K_0$ is split.  The inclusion of $F^5 = \Z$ in $K_0$ is split by the Euler characteristic $\chi\colon K_0(X) \to \Z$.  The inclusion of $F^4/F^5 = \Z$ in $K_0/F^5$ is split by $\chi(\underline\quad \otimes \O_H)$, where $H \subset X$ is a hyperplane section, because $CH^4$ is generated by a line.  The rest of the filtration splits because the quotients are all $\Z$.  Thus
\[ K_0(X) = \Z^6 \oplus IJ(X). \]

The remainder of the proof proceeds like the proofs of Proposition \ref{3fold_KAX} and Theorem \ref{4fold_KAX}(b).
\end{proof}

Projecting from a plane in the cubic 5-fold $X$ yields a quadric surface bundle over $\P^3$ that degenerates over a sextic surface with 31 nodes.  It would be interesting to study the big resolution $\tilde Y$ of the double cover of $\P^3$ branched over this sextic, which should carry a Brauer class $\alpha$ of order 2; we suspect that $\dimQL(\tilde Y) = 2$ but $\dimQL(\tilde Y,\alpha) = 1$.